\documentclass[journal]{new-aiaa}
\usepackage[utf8]{inputenc}

\usepackage{graphicx}
\usepackage{amsmath}
\usepackage[version=4]{mhchem}
\usepackage{siunitx}
\usepackage{longtable,tabularx}
\usepackage{indentfirst}
\usepackage{hyperref}
\usepackage{mathtools}
\usepackage{bm}
\usepackage{xspace}
\usepackage{wrapfig}
\usepackage{multicol}
\usepackage{xurl}
\usepackage[flushleft]{threeparttable}
\usepackage{booktabs}
\usepackage{gensymb}
\usepackage{mdframed}
\usepackage{subcaption}
\usepackage{xcolor, pifont}
\usepackage{lscape}
\usepackage{tcolorbox}
\usepackage[ruled, vlined, linesnumbered]{algorithm2e}

\let\svthefootnote\thefootnote
\newcommand\freefootnote[1]{
  \let\thefootnote\relax
  \footnotetext{#1}
  \let\thefootnote\svthefootnote
}

\newcommand{\MCRP}{\hyperlink{MCRP}{\textsf{MCRP}}\xspace}
\newcommand{\SMCRP}{\hyperlink{SMCRP}{\textsf{SMCRP}}\xspace}
\newcommand{\SDDiP}{\hyperlink{SDDiP}{\textsf{SDDiP}}\xspace}

\definecolor{myblue}{rgb}{0, 0.23, 0.64}
\definecolor{WVUblue}{rgb}{0, 0.16, 0.33}
  
\hypersetup{
    colorlinks=true,
    linkcolor=myblue,
    filecolor=magenta,
    citecolor=myblue,
    urlcolor=myblue,
}

\title{Stochastic Multistage Constellation Reconfiguration Problem: Stochastic Dual Dynamic Integer Programming Approach}

\author{Brycen D. Pearl\footnote{Ph.D. Candidate, Department of Mechanical, Materials and Aerospace Engineering.} and Hang Woon Lee\footnote{Assistant Professor, Department of Mechanical, Materials and Aerospace Engineering; hangwoon.lee@mail.wvu.edu. Member AIAA (Corresponding Author).}}
\affil{West Virginia University, Morgantown, WV, 26506}

\begin{document}

\newpage

\maketitle

\begin{abstract}
    Observation tasking is critical to satellite operations, enabling observation of planetary phenomena, orbital debris monitoring, and space domain awareness. To improve the effectiveness of satellite operations, orbital maneuverability is introduced as a leading-edge concept of operations within constellation reconfigurability for response to dynamic events. Previous investigations regarding constellation reconfigurability consist of deterministic mission environments, in which the formulations consider \textit{a priori} knowledge of the environment; however, observation objectives inherently involve uncertainty at any given time. As such, scheduling satellite operations must account for uncertainties to ensure adequate observation tasking. In response, we present a stochastic variant of the Multistage Constellation Reconfiguration Problem (MCRP) which is solved using Stochastic Dual Dynamic Integer Programming (SDDiP). We additionally explore other stochastic problem-solving techniques for solution quality comparison. To demonstrate each solution method, two computational experiments with stochastic target properties are conducted. The first concerns random orbital targets, and the second concerns simulated hurricanes. The results of the experiments demonstrate the effectiveness of the SDDiP solution approach over other stochastic problem-solving methods. Overall, the stochastic MCRP accounts for target stochasticity while obeying visible time windows and maneuver feasibility.
\end{abstract}

\section*{Nomenclature}
{
\renewcommand\arraystretch{1.0}
\noindent\begin{longtable*}{@{}l @{\quad=\quad} l@{}}
$\mathcal{A}$        & Markov decision process action space, index $a$, cardinality $J^K$ \\
$c$                  & cost of orbital maneuver \\
$c_{\max}$           & budget for total orbital maneuver costs \\
$E$                  & probability of scenario \\
$I$                  & \SDDiP iteration number \\
$\mathcal{J}$        & set of orbital slots, indices $i,j$, cardinality $J$ \\
$\mathcal{K}$        & set of satellites, index $k$, cardinality $K$ \\
$L$                  & lower bound of convex lower approximation \\
$\mathcal{L}$        & \SMCRP Lagrangian relaxation \\
$M$                  & number of sampled scenarios \\
$\mathcal{N}$        & Markov decision process state space, indices $n, n^\prime$, cardinality $1+S(J^K)W$ \\
$P_I$                & \SMCRP forward problem \\
$\mathcal{P}$        & set of targets, index $p$, cardinality $P$ \\
$\mathcal{Q}$        & Markov decision process transition function \\
$R_I$                & suitable \SMCRP relaxation \\ 
$\mathcal{R}$        & Markov decision process reward function \\
$\mathcal{S}$        & set of reconfiguration stages, index $s$, cardinality $S+1$ \\
$T_r$                & finite schedule duration \\
$\mathcal{T}$        & set of time steps, index $t$, cardinality $T$ \\
$V$                  & binary visibility condition of targets \\
$\bar{v}$            & Integer Optimality Cut objective values \\
$\tilde{v}$          & Strengthened Benders' Cut objective values \\
$\bar{x}$            & Integer Optimality Cut solution location \\
$\mathcal{W}$        & set of scenarios, index $w$, cardinality $W$ \\
$x$                  & decision variable to control constellation reconfiguration \\
$y$                  & indicator variable of obtained target visibility \\
$Z$                  & sum of rewards from given policy \\
$z$                  & auxiliary stage linking variable \\
$z_{\alpha/2}$       & standard normal critical value \\
$\Delta t$           & time step size between discrete time steps \\
$\pi$                & observation reward of target \\
$\theta$             & auxiliary lower approximation variable \\
$\tilde{\lambda}$    & Strengthened Benders' Cut Lagrangian dual multiplier values \\ 
$\xi$                & realization of uncertain parameters \\
$\psi$               & convex lower approximation of function value \\
$\Omega$             & collection of sampled uncertain parameter values \\
\end{longtable*}
}

\section{Introduction}

Satellite constellations are often required to observe targets whose locations, priorities, and values evolve over the course of a mission. This requirement arises in Earth-observation and space-domain-awareness applications, including natural-disaster monitoring \cite{Panteras2018Flood,Kaku2019Disaster} and the detection and tracking of orbital debris or spacecraft \cite{Zhang2024Debris,Colombi2018SSA}. Because orbital geometry constrains when and where each satellite can collect observations, constellation performance depends not only on scheduling available observation opportunities but also on adapting the constellation geometry to changing mission needs. This challenge becomes particularly important when target behavior and environmental conditions are uncertain, since maneuvering decisions must be made sequentially before all future information is known.

\textit{Constellation reconfigurability} provides a mechanism for such adaptation by enabling satellites to perform coordinated orbital maneuvers that reshape the constellation in support of a specified mission objective \cite{Siddiqi2005Optimal}. Prior research has shown that reconfigurable constellations can outperform both fixed-orbit, non-maneuverable constellations and constellations of agile satellites limited to attitude control \cite{Pearl2025Benchmarking}. Beyond these comparative studies, constellation reconfigurability has been investigated in several distinct mission-planning problems, including minimizing the maximum revisit time for emergency targets using sequential quadratic programming \cite{Fang2026Reconfiguration} and planning collision-aware orbital maneuvers through a Markov decision process formulation \cite{Xu2024RLReconfiguration}. Of particular importance is multistage constellation reconfiguration, the framework for modeling and optimizing the sequence of orbital maneuvers through multiple opportunities, or stages, during the given mission duration \cite{Lee2022Maximizing}. An implementation of constellation reconfigurability in the environment of the Earth Observation Satellite Scheduling Problem (EOSSP) is detailed in Ref.~\cite{Pearl2026REOSSP}, which schedules the satellite tasks of orbital maneuvers, target observation, data downlink to ground stations, and solar charging, while tracking onboard data and battery storage. To the best of the author's knowledge, all investigations into the concept of constellation reconfigurability have assumed a \textit{deterministic} operational setting, in which knowledge of the environment is known before the start of operations, apart from one investigation detailed in Ref.~\cite{Pearl2025Stochastic}.

Operating satellite systems requires consideration of uncertainty, as many natural and human-made ground- and space-based phenomena rely on complex processes that are difficult to predict in advance. For instance, natural disasters are complex systems governed by uncertain dynamics, such as a volcano that may erupt at times that are difficult to determine \cite{Chien2020Automated}, or a hurricane that may change in intensity or heading by large amounts \cite{Ruf2019CYGNSS}. Uncertainty also affects space domain awareness involving orbital objects such as satellites and debris. Many debris objects are too small to maintain adequate tracking custody \cite{Liou2011Debris}, and errors in orbital propagation may accumulate due to uncertainties in atmospheric drag and solar radiation pressure models \cite{Vallado2014Drag}. As such, satellite task scheduling must account for the inherent uncertainty of important targets to ensure reliable operations, an aspect that remains largely unaddressed in constellation-reconfigurability research.

To address uncertainty in the decision-making process for satellite tasking, stochastic decision-making methods are leveraged to ensure that effective decisions are made regarding probabilistic outcomes \cite{Kochenderfer2015SPS}. One such mathematical framework is the Markov decision process (MDP), a model for sequential decision-making under uncertainty relying on the Markovian nature that the outcome of a decision is independent of previous decisions once the current state is known \cite{Puterman1994MDPVI}. For instance, Ref.~\cite{Nastasi2019POMDPtasking} leverages a partially observable MDP (POMDP) to control the tasking of ground- and space-based sensors to track maneuvering spacecraft using various sensor types with uncertainty in the maneuvering target. Similarly, Ref.~\cite{Fedeler2022Cislunar} utilizes a POMDP to task the pointing of cislunar satellites towards various cislunar targets with uncertainty derived from the individual observations obtained throughout a mission. MDP formulations have also been leveraged to obtain schedule solutions from the Agile EOSSP, a scheduling problem where satellites may perform attitude control to gain observations of ground targets \cite{Wang2024Strategy}. Key research in this area includes the use of a semi-MDP solved via deep reinforcement learning \cite{Stephenson2025POMDP} and Deep Q-Learning neural network techniques to solve a POMDP formulation \cite{He2022DQL}.

To integrate constellation reconfigurability while accounting for the uncertainty inherent in the operation of satellite systems, this paper develops the \textit{Stochastic Multistage Constellation Reconfiguration problem} (SMCRP) by building upon the author's previous research on the Multistage Constellation Reconfiguration Problem (MCRP) \cite{Lee2022Maximizing,Lee2023NovelFormulation,Lee2024Deterministic}. The MCRP is a mixed integer linear programming (MILP) problem that maximizes the obtained target observation rewards from a constellation of maneuverable satellites with respect to a deterministic operation environment. The MCRP accounts for target visibility, orbital maneuver feasibility, and maneuver budget constraints while optimally selecting orbital maneuvers without uncertainty. To overcome the deterministic assumptions of the MCRP, the SMCRP incorporates various scenarios wherein many operational parameters contain uncertain values governed by probabilistic realization in each scenario, thus implementing stochasticity under similar orbital maneuvering constraints. The objective of the SMCRP is to maximize the expected observation rewards under environmental uncertainty through the use of Stochastic Dual Dynamic Integer Programming (SDDiP), returning an optimal policy of orbital maneuvers to be performed sequentially relative to the realization of uncertain parameters \cite{Zou2019SDDiP}. The concept of SDDiP extends upon Stochastic Dual Dynamic Programming, an optimization technique that decomposes larger sequential decision-making problems into a series of smaller subproblems that are solved recursively in iterations \cite{Yuan2021StochasticDP}, to include integer constraints on a subset of decision variables \cite{Zou2019SDDiP}. 

To validate the use of SDDiP to solve the SMCRP, this paper additionally conducts computational experiments relative to orbital targets that are deployed into random orbits at random times during a mission, as well as simulated hurricanes with probabilistic travel trajectories. The computational experiments provide a comparison between the SDDiP approach and other approaches, such as the use of random maneuvers, a non-maneuverable constellation, solving an MDP reformulation, and a deterministic upper bound set by the MCRP. The results indicate that the SDDiP approach outperforms all other solution methods, approaching the deterministic upper bound in many cases. The SMCRP obeys principles of target visibility and orbital maneuver feasibility, while accounting for the stochastic nature of the mission environment. 

\begin{figure}[!ht]
    \centering
    \includegraphics[width=0.75\textwidth]{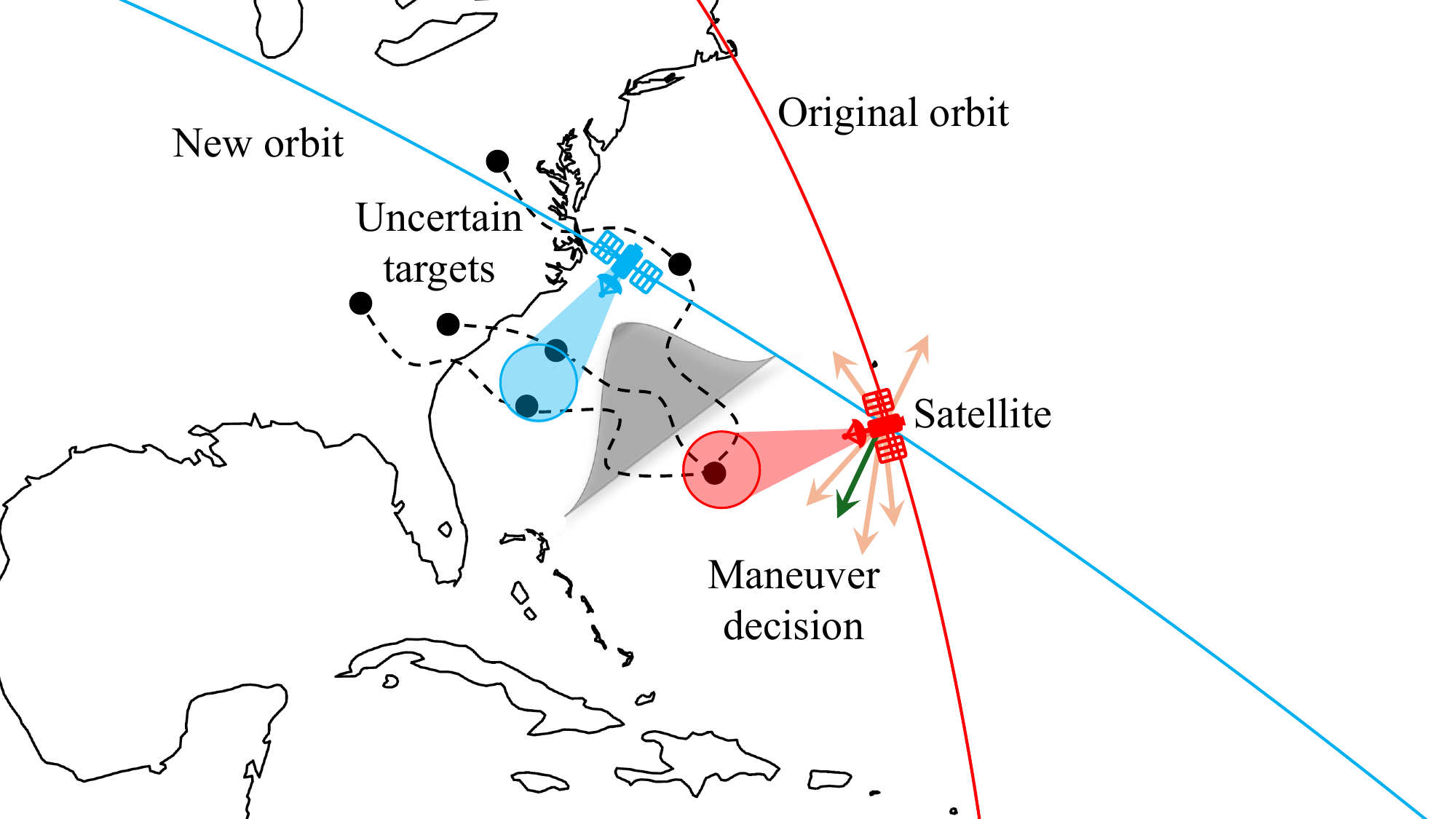}
    \caption{A demonstration of the SMCRP.}
    \label{fig:Overarching_Stochastic}
\end{figure}

An illustration of the purpose and result of the developed SMCRP is shown in Fig.~\ref{fig:Overarching_Stochastic}. Within Fig.~\ref{fig:Overarching_Stochastic}, an uncertain ground target is shown progressing from an initial position along various trajectories over time with a denoted probability density function representing the probability of each trajectory. The SMCRP leverages the uncertain target points over time with their associated probability to determine the best orbital maneuver for a satellite or constellation of satellites that maximizes the obtained visibility of the uncertain target. Figure~\ref{fig:Overarching_Stochastic} demonstrates the various orbital maneuver options available to the satellite, wherein the obtained visibility of the uncertain target is maximized after the maneuver is performed. 

The rest of the paper is organized as follows. First, Sec.~\ref{sec:DMCRP} reviews the \MCRP from Refs.~\cite{Lee2022Maximizing,Lee2023NovelFormulation,Lee2024Deterministic} with deterministic target setting, and Sec.~\ref{sec:SMCRP} details the modifications made to the \MCRP for the incorporation of stochastic target scenarios alongside the incorporation of SDDiP as a solution method. Then, Sec.~\ref{sec:Experimentation} details additional stochastic decision-making approaches and two computational experiment environments, the first regarding monitoring orbital targets and the second concerning simulated hurricanes. Finally, Sec.~\ref{sec:Conclusion} discusses the overall conclusions provided through the case study results and offers potential avenues for future research. 

\section{Deterministic Multistage Constellation Reconfiguration Problem} \label{sec:DMCRP}

Previous research provided in Refs.~\cite{Lee2022Maximizing,Lee2023NovelFormulation,Lee2024Deterministic} depicts the development of the \MCRP, wherein Ref.~\cite{Lee2024Deterministic} contains the most recent version, which is expanded upon in this paper. 

The \MCRP maximizes the obtained observation rewards of a given set of targets gained by a constellation of satellites capable of performing orbital maneuvers for constellation reconfiguration over equally spaced stages within a given mission horizon. The mission horizon is represented as a set of discrete time steps $t \in \mathcal{T} = \{1, 2, \ldots, T\}$ with $T = T_r/\Delta t$ total time steps, where $T_r \ge 0$ is the finite mission duration and $\Delta t > 0$ is the time step size. Within the mission horizon, $S > 0$ total stages are equally spaced at intervals of $T^s = T/S$ and the stage horizon is defined as $\mathcal{T}^s$ such that $\mathcal{T} = \{\mathcal{T}^1, \mathcal{T}^2, \ldots, \mathcal{T}^S\}$. It should be noted that the initial constellation configuration is defined at stage $s = 0$, which occurs before the start of the mission horizon. As such, the set of stages is defined as $\mathcal{S} = \{0, 1, 2, \ldots, S\}$, where each stage $s \in \mathcal{S} \setminus \{0\}$ is associated with the time steps $\mathcal{T}^s = \{t+(s-1)T^s \colon t = 1, 2, \ldots, T^s\}$. Additionally, the set of satellites within the constellation is defined as $\mathcal{K} = \{1, 2, \ldots, K\}$ with $K > 0$ total satellites and associated orbital slots $\mathcal{J}^{sk} = \{1, 2, \ldots, J^{sk}\}, s \in \mathcal{S}, k \in \mathcal{K}$ where $J^{sk} \ge 0$ is the number of orbital slots. While the \MCRP in Refs.~\cite{Lee2022Maximizing,Lee2023NovelFormulation,Lee2024Deterministic} considers a heterogeneous set of satellites, in which only a subset are capable of performing maneuvers, this work assumes a homogeneous set of satellites such that all are capable of performing maneuvers. The final set is that of targets, defined as $\mathcal{P} = \{1, 2, \ldots, P\}$ with $P \ge 0$ total targets and associated locations and time-dependent observation rewards $\pi^s_{tp} \in \mathbb{R}_+$. 

The \MCRP returns two optimal variables: a binary variable that represents the sequence of orbital maneuvers for each satellite over all stages within the mission horizon, and a continuous variable that represents the obtained visibility of each target by the constellation as a whole. The binary decision variable is $x^{sk}_{ij}$, which indicates the orbital maneuvers of satellites; $x^{sk}_{ij} = 1$ if satellite $k \in \mathcal{K}$ performs an orbital maneuver from orbital slot $i \in \mathcal{J}^{s-1, k}$ in the previous stage $s-1$ to orbital slot $j \in \mathcal{J}^{sk}$ in the current stage $s$ for stage $s \in \mathcal{S} \setminus \{0\}$ and $x^{sk}_{ij} = 0$ otherwise. Each orbital maneuver has an associated cost, denoted $c^{sk}_{ij} \ge 0$ with indexing identical to $x^{sk}_{ij}$, and each satellite $k$ has an associated maximum maneuver budget allotted for the entire mission horizon, denoted $c^k_{\max} \ge 0$. Separately, the continuous indicator variable is $y^s_{tp} \in [0, 1]$, which indicates the obtained visibility of target $p \in \mathcal{P}$ at time $t \in \mathcal{T}^s$ of stage $s \in \mathcal{S} \setminus \{0\}$. Crucially, Refs.~\cite{Lee2022Maximizing,Lee2023NovelFormulation,Lee2024Deterministic} consider some positive integer number of satellites as a coverage requirement to obtain target visibility, while earlier work in Ref.~\cite{Lee2023Regional} clarifies that $y^s_{tp}$ may be continuous on the range from zero to one if said coverage requirement is limited to one. As such, the \MCRP in this work assumes that only one satellite in the constellation is required to have visibility of a target to contribute to the obtained visibility $y^s_{tp}$, subject to the binary VTW, denoted $V^{sk}_{tjp}$; $V^{sk}_{tjp} = 1$ if satellite $k$, located in orbital slot $j \in \mathcal{J}^{sk}$ during stage $s \in \mathcal{S} \setminus \{0\}$ at time $t \in \mathcal{T}^s$, has some sensor access to target $p$, and $V^{sk}_{tjp} = 0$ otherwise. 

The mathematical formulation of the \MCRP is shown in Formulation~\eqref{MCRP} \cite{Lee2022Maximizing,Lee2023NovelFormulation,Lee2024Deterministic}: \hypertarget{MCRP}{}
\begin{subequations}
    \begin{alignat}{2}
\max \quad & \sum_{s \in \mathcal{S} \setminus \{0\}} \sum_{t \in \mathcal{T}^s} \sum_{p \in \mathcal{P}} \pi^s_{tp} y^s_{tp}
\label{MCRP:obj} \\
\text{s.t.} \quad & \sum_{j \in \mathcal{J}^{1k}} x^{1k}_{ij} = 1, \qquad && \forall k \in \mathcal{K}, \forall i \in \mathcal{J}^{0k}
\label{MCRP:init} \\
& \sum_{j \in \mathcal{J}^{s+1, k}} x^{s+1, k}_{ij} - \sum_{j^\prime \in \mathcal{J}^{s-1}} x^{sk}_{j^\prime i} = 0 , \qquad && \forall s \in \mathcal{S} \setminus \{0, S\}, \forall k \in \mathcal{K}, \forall i \in \mathcal{J}^{sk}
\label{MCRP:flow} \\
& \sum_{s \in \mathcal{S} \setminus \{0\}} \sum_{i \in \mathcal{J}^{s-1, k}} \sum_{j \in \mathcal{J}^{sk}} c^{sk}_{ij} x^{sk}_{ij} \le c^k_{\max}, \qquad && \forall k \in \mathcal{K} 
\label{MCRP:cost} \\
& \sum_{k \in \mathcal{K}} \sum_{i \in \mathcal{J}^{s-1, k}} \sum_{j \in \mathcal{J}^{sk}} V^{sk}_{tjp}x^{sk}_{ij} \ge y^s_{tp}, \qquad && \forall s \in \mathcal{S} \setminus \{0\}, \forall t \in \mathcal{T}^s, \forall p \in \mathcal{P}
\label{MCRP:vis} \\
& x^{sk}_{ij} \in \{0,1\}, && \forall s \in \mathcal{S} \setminus \{0\}, \forall k \in \mathcal{K}, \forall i \in \mathcal{J}^{s-1, k}, \forall j \in \mathcal{J}^{sk} 
\label{MCRP:x} \\
& y^s_{tp} \in [0,1], && \forall s \in \mathcal{S} \setminus \{0\}, \forall t \in \mathcal{T}^s, \forall p \in \mathcal{P} 
\label{MCRP:y}
    \end{alignat}
    \label{MCRP}
\end{subequations}

Objective function \eqref{MCRP:obj} maximizes the obtained observation reward through an optimal sequence of orbital maneuvers. Constraints~\eqref{MCRP:init} control the orbital maneuvers of each satellite from the initial conditions $\mathcal{J}^{0k}$ such that only one arrival orbital slot $j \in \mathcal{J}^{1k}$ is selected in the first stage. Similarly, constraints~\eqref{MCRP:flow} control all subsequent orbital maneuvers for stages $s \in \mathcal{S} \setminus \{0, S\}$ such that any orbital maneuver from orbital slot $i \in \mathcal{J}^{sk}$ to orbital slot $j \in \mathcal{J}^{s+1, k}$ may only occur if the previous orbital maneuver from orbital slot $j^\prime \in \mathcal{J}^{s-1, k}$ resulted in arrival to orbital slot $i$. Additionally, constraints~\eqref{MCRP:cost} apply the cost $c^{sk}_{ij}$ to each orbital maneuver such that the associated budget $c^k_{\max}$ is not exceeded. Separately, constraints~\eqref{MCRP:vis} apply the VTW condition to $y^s_{tp}$, ensuring visibility according to the currently occupied orbital slot $j \in \mathcal{J}^{sk}$. Finally, constraints~\eqref{MCRP:x} and~\eqref{MCRP:y} define the binary and continuous domains of the decision and indicator variable, respectively.

\section{Stochastic Multistage Constellation Reconfiguration Problem} \label{sec:SMCRP}

Expanding upon the \MCRP, the SMCRP incorporates stochastic target properties through the use of discrete target scenarios, leveraging SDDiP to obtain an optimal orbital maneuver policy in each scenario through concepts of stochastic dynamic programming, mixed integer programming, and the overall problem structure. 

The SMCRP maximizes the expected observation rewards obtained from a given set of probabilistic targets using constellation reconfigurability in the same manner as the \MCRP, returning an effective orbital maneuver sequence. As such, the parameters of the SMCRP are similar in nature to those of the \MCRP, with adaptation to reflect the uncertain nature of the target set. Specifically, the SMCRP incorporates a set of $W$ total scenarios, $\mathcal{W} = \{1, 2, \ldots, W\}$, with unique target set $\mathcal{P}_w = \{1, \ldots, P_w\}$. Each scenario $w\in\mathcal{W}$ has probability $E_w$ defined by some probability mass function such that $\sum_{w\in\mathcal{W}}E_w = 1$. As a result of the additional index for the scenario of each target set, the VTW parameter is expanded to $V^{sk}_{tjpw} \in \{0, 1\}$ and the time-dependent observation reward parameter is expanded to $\pi^s_{tpw} \in \mathbb{R}_+$, wherein the target $p \in \mathcal{P}_w$ refers to a target in scenario $w \in \mathcal{W}$. The SMCRP may be reformulated to fit the structure and assumptions necessary for SDDiP, for which SDDiP itself and any necessary modifications to the SMCRP are subsequently described.

\subsection{Stochastic Dual Dynamic Integer Programming} \label{subsec:SDDiP}

The concept of SDDiP, as described in Ref.~\cite{Zou2019SDDiP}, utilizes the theory of stochastic dual dynamic programming for solving multistage stochastic programs alongside key adaptations to fit problems that contain a subset of integer variables. A multistage stochastic integer program (MSIP) is defined as follows, using the native notation of Ref.~\cite{Zou2019SDDiP} disjoint from the notation throughout the rest of this paper: 
\begin{subequations}
    \begin{alignat}{2}
        \min \quad & \sum_{n \in \mathcal{T}} p_nf_n(x_n, y_n) \\
        \text{s.t.} \quad & (z_n, x_n, y_n) \in X_n, \quad && \forall n \in \mathcal{T} \\
        & z_n = x_{a(n)}, && \forall n \in \mathcal{T} \label{MSIP:Linking} \\
        & z_n \in \left[ 0, 1 \right]^d, && \forall n \in \mathcal{T} \\
        & x_n \in \{0, 1\}^d, && \forall n \in \mathcal{T}
    \end{alignat}
    \label{MSIP}
\end{subequations}

Within formulation~\eqref{MSIP}, $\mathcal{T}$ is the scenario tree with nodes denoted by $n$, each of which is linked to parent node $a(n)$ by the binary state variable $x_n$ with dimension $d$. Moreover, $y_n$ is the local variable within node $n$, $z_n$ is the local copy of the parent node state, and $p_n$ is the probability of node $n$ occurring in the scenario tree. Formulation~\eqref{MSIP} is linear, such that $f_n(x_n, y_n)$ is the linear objective function, and $X_n$ is the constraint set of the form $B_nz_n + A_nx_n + C_ny_n \ge b_n$ with integer constraints on a subset of $y_n$. While constraints~\eqref{MSIP:Linking} may seem redundant initially, they are crucial to later aspects of the SDDiP algorithm. The solution of formulation~\eqref{MSIP} returns a policy, that being the binary state decision variables, at each node of the scenario tree, which is a direct mapping of a decision to take in each possible realization of scenarios. The MSIP in formulation~\eqref{MSIP} is typically large in scale, as the scenario tree grows exponentially with the dimension of the uncertain parameters and decision-making stages; as such, formulation~\eqref{MSIP} is reformulated via the following dynamic programming equations: 
\begin{subequations}
    \begin{alignat}{2}
        P_1 \coloneq \quad & \min \quad &&f_1(x_1, y_1) + \sum_{m \in \mathcal{C}(1)} q_{1m}Q_m(x_1) \\
        & ~\text{s.t.} \quad && (z_1, x_1, y_1) \in X_1 \\
        & \quad && z_1 = x_{a(1)},~ z_1 \in\left[0, 1\right]^d,~ x_1 \in\{0, 1\}^d\\
        P_n, ~ \forall n \in \mathcal{T} \setminus \{1\} \coloneq \quad  Q_n(x_{a(n)}) = &\min  \quad &&f_n(x_n, y_n) + \sum_{m \in \mathcal{C}(n)} q_{nm}Q_m(x_n) \\
        & ~\text{s.t.} \quad && (z_n, x_n, y_n) \in X_n \\
        & \quad && z_n = x_{a(n)},~ z_n \in\left[0, 1\right]^d,~ x_n \in\{0, 1\}^d
    \end{alignat}
\end{subequations}
where $q_{nm} = p_m/p_n$ is the conditional probability of transition from node $n$ to it's child $m \in \mathcal{C}(n)$. Additionally, $Q_n(\cdot)$ is the value function at node $n$, representing the impact of future decisions within each child node problem. For the conditions of SDDiP to hold, the following assumptions are made regarding the MSIP: 1) the first stage and the initial conditions $x_0$ are deterministic, 2) the objective functions $f_n(x_n, y_n)$ are linear, 3) the constraint sets $X_n$ are nonempty compact mixed integer polyhedral sets, 4) the state variables $x_n$ are binary, 5) the MSIP has complete continuous recourse, and 6) the MSIP has stagewise independence. By complete continuous recourse, there exists a feasible next decision $(x_n, y_n)$ given any previous decision $x_{a(n)}$, and by stagewise independence, the uncertain parameters only depend on the decision-making stage rather than the node itself (thus being independent of previous decisions). Regarding the final assumption on stagewise independence, Refs.~\cite{Zou2019SDDiP,Zou2019UnitCommitment} define a reformulation of the dynamic programming equations for use in iteration $i$ of the SDDiP algorithm as follows: 
\begin{subequations}
    \begin{alignat}{2}
        P^i_t \coloneq \quad  Q^i_t(x^i_{t-1}, \psi^i_t, \xi^k_t) = &\min  \quad &&f_t(x_t, y_t) + \psi^i_t(x_t) \\
        & ~\text{s.t.} \quad && (z_t, x_t, y_t) \in X_t \\
        & \quad && z_t = x^i_{t-1},~ z_t \in\left[0, 1\right]^d,~ x_t \in\{0, 1\}^d \\
        \psi^i_t \coloneq &\min\quad &&  \theta_t \\
        & ~\text{s.t.}&& \theta_t \ge L_t \\
        & && \theta_t \ge \text{cuts}^i \label{Cut_eq}
    \end{alignat}
\end{subequations}
where $\psi_t$ is the convex lower approximation of the value function at decision $t$, $\xi^k_t$ is the $k$th realization of the uncertain parameters at decision $t$, and constraints~\eqref{Cut_eq} refer to the various cuts to improve the lower approximation in iterations of the SDDiP algorithm. The cuts in constraints~\eqref{Cut_eq} are described in more detail later in this subsection.

The \MCRP may be reformulated to fit the MSIP dynamic programming equations, with the notation native to Sec.~\ref{sec:DMCRP}, through slight modification to appropriately fit the aforementioned assumptions. The reformulated dynamic programming equations defined in the \SMCRP are as follows: \hypertarget{SMCRP}{}
\begin{subequations} 
    \begin{alignat}{2}
\begin{split}
    P^s_{I}\left( x^{s-1,k}_{ij;Iw}, \psi^s_{I}, \xi^s_w \right)&, \forall s \in \mathcal{S} \setminus\{0\} \coloneq\\
    & \min \qquad \left(-\sum_{t \in \mathcal{T}^s} \sum_{p \in \mathcal{P}_w} \pi^s_{tpw}y^s_{tp}\right) + \psi^s_{I}\left( x^{sk}_{ij} \right) 
\end{split}
\label{SMCRP:obj}\\
&\text{s.t.} \quad \sum_{j \in \mathcal{J}^{sk}} x^{sk}_{ij} - \sum_{j^\prime \in \mathcal{J}^{s-2, k}} z^{sk}_{j^\prime i} = 0, \quad && \forall k \in \mathcal{K}, \forall i \in \mathcal{J}^{s-1, k} 
\label{SMCRP:Flow}\\
& \qquad x^{s-1,k}_{ij;Iw} = z^{sk}_{ij}, \quad && \forall k \in \mathcal{K}, \forall i \in \mathcal{J}^{s-2, k}, \forall j \in \mathcal{J}^{s-1, k}
\label{SMCRP:Link}\\
& \qquad \sum_{k \in \mathcal{K}} \sum_{i \in \mathcal{J}^{s-1,k}} \sum_{j \in \mathcal{J}^{sk}} V^{sk}_{tjpw} x^{sk}_{ij} \ge y^s_{tp}, \quad && \forall t \in \mathcal{T}^s, \forall p \in \mathcal{P}_w
\label{SMCRP:vis}\\
& \qquad x^{sk}_{ij} \in \{0,1\}, \quad && \forall k \in \mathcal{K}, \forall i \in \mathcal{J}^{s-1, k}, \forall j \in \mathcal{J}^{sk} 
\label{SMCRP:x} \\
& \qquad z^{sk}_{ij} \in [0,1], \quad && \forall k \in \mathcal{K}, \forall i \in \mathcal{J}^{s-2, k}, \forall j \in \mathcal{J}^{s-1,k} 
\label{SMCRP:z} \\
& \qquad y^s_{tp} \in [0,1], \quad && \forall t \in \mathcal{T}^s, \forall p \in \mathcal{P}_w 
\label{SMCRP:y}
    \end{alignat}
    \label{SMCRP}
\end{subequations}

Within the \SMCRP, many parameters and variables inherit their purpose from the \MCRP; note that some contain a semicolon (;) in the subscript. The indices to the left of the semicolon refer to dimensions of the parameter or variable, while indices to the right are used to convey algorithm information. Firstly, the binary state variable linking successive decision-making stages is the decision variable $x^{sk}_{ij}$ that controls the orbital maneuver sequence. The binary state variable links each stage through constraints~\eqref{SMCRP:Flow} and~\eqref{SMCRP:Link}; the former constraints ensure the logical continuity of orbital maneuvers in the same manner as constraints~\eqref{MCRP:init} and~\eqref{MCRP:flow}, and the latter assign the auxiliary local copy $z^{sk}_{ij}$ of the previous stage decision $x^{s-1,k}_{ij;Iw}$ dependent upon iteration $I$ of the algorithm and scenario $w \in \mathcal{W}$ of the uncertain parameters. Second, the local variable is the indicator variable $y^s_{tp}$ that determines target visibility on the continuous range from zero to one. The local variable indicates the visibility of targets via constraints~\eqref{SMCRP:vis} in the same manner as constraints~\eqref{MCRP:vis}, under the consideration that only one satellite is required to obtain visibility. Constraints~\eqref{SMCRP:x}--\eqref{SMCRP:y} define the domain of the stage and local variables in each stage subproblem $P^s_{I}\left( x^{s-1,k}_{ij;Iw}, \psi^s_{I}, \xi^s_w \right)$. Finally, $\xi^s_w$ is the realization of $\pi^s_{tpw}$ and $V^{sk}_{tjpw}$ in stage $s$ and scenario $w$, and $\psi^s_{I}\left(\cdot\right)$ is defined: 
\begin{subequations}
    \begin{alignat}{2}
\psi^s_{I}\left( x^{sk}_{ij} \right)  \coloneq & \quad \min \quad && \theta^s 
\label{SMCRP_psi:obj}\\
& \quad \text{s.t.} \quad && \theta^s \ge L^s 
\label{SMCRP_psi:LB}\\
& && \theta^s \ge \text{cuts}_{I}
\label{SMCRP_psi:cuts}
    \end{alignat}
    \label{SMCRP_psi}
\end{subequations}
where $L^s$ in constraint~\eqref{SMCRP_psi:LB} is the lower bound of the lower approximation $\psi^s_{I}\left(\cdot\right)$, and constraints~\eqref{SMCRP_psi:cuts} improve the lower approximation by providing adequate cut information in iteration $I$. It should be noted that the constraints~\eqref{MCRP:cost} are not included in the \SMCRP due to Assumption 6 on stagewise independence. Instead, a consideration is made that the orbital slots provided to each satellite in $\mathcal{J}^{sk}, \forall s \in \mathcal{S} \setminus \{0\}, \forall k \in \mathcal{K}$ are constructed such that the costs incurred in $c^{sk}_{ij}$ by any decisions $x^{sk}_{ij}$ do not exceed a desired budget. A more detailed explanation of such a construction of the orbital slots is located in Appendix A. Otherwise, all other assumptions are met inherently: 1) the initial conditions $x^{0k}_{ij}$ may be set as deterministic using a slight abuse of notation; $x^{0k}_{ij} \equiv j \in \mathcal{J}^{0k}$ where
\begin{equation*}
    x^{0k}_{ij} = \begin{cases}
        1, & \text{if } i = j \in \mathcal{J}^{0k} \\
        0, & \text{otherwise}
    \end{cases}, \quad \forall k \in \mathcal{K}
\end{equation*}
2) the objective function~\eqref{SMCRP:obj} is linear, 3) constraints~\eqref{SMCRP:Flow}--\eqref{SMCRP:y} are nonempty ($K > 0$), compact (all variables are bounded), mixed integer (binary nature of $x^{sk}_{ij}$), and polyhedral (each constraint is a linear inequality or equality), 4) the state variables $x^{sk}_{ij}$ are binary, and 5) the \SMCRP has complete continuous recourse. The justification of Assumption 5 warrants further explanation; the \SMCRP has a feasible solution given any previous decision $x^{s-1,k}_{ij;Iw}$ via the following. First, via constraints~\eqref{SMCRP:Link}, there is a feasible local copy of the previous decision $z^{sk}_{ij} \in [0, 1]$ since $x^{s-1,k}_{ij;Iw}$ is binary. Second, via constraints~\eqref{SMCRP:Flow}, there is a feasible current decision $x^{sk}_{ij} = 1, \text{ for any } j \in \mathcal{J}^{sk}$ where $i\in\mathcal{J}^{s-1,k}$ satisfies $z^{sk}_{j^\prime i}=1$ relative to each satellite $k\in\mathcal{K}$. Finally, via constraints~\eqref{SMCRP:vis}, $y^s_{tp} = 0 \in [0, 1], ~\forall t \in \mathcal{T}^s, ~\forall p \in \mathcal{P}_w$ is guaranteed to be feasible, as the quantity $\sum_{k \in \mathcal{K}} \sum_{i \in \mathcal{J}^{s-1,k}} \sum_{j \in \mathcal{J}^{sk}} V^{sk}_{tjpw} x^{sk}_{ij}$ is binary. As such, the \SMCRP meets all conditions of SDDiP and may be applied to solve the \SMCRP and return an optimal set of maneuvers as a policy for each scenario.

The \SDDiP algorithm, originally depicted in Ref.~\cite{Zou2019SDDiP}, is restated in Algorithm~\ref{alg:SDDiP} using the notation of the \SMCRP. In each iteration $I$ of the \SDDiP algorithm, five steps occur: the Sampling step, Forward step, Upper bound update, Backward step, and Lower bound update. Before iterations occur, the \SDDiP algorithm initializes the upper bound $UB$, lower bound $LB$, iteration count $I$, and lower approximation $\psi^s_1\left(\cdot\right) \equiv 0,~\forall s \in \mathcal{S}\setminus\{0\}$ for use in the first iteration. The \SDDiP algorithm continues iterating until some convergence criteria are satisfied, where Ref.~\cite{Zou2019SDDiP} suggests stopping when the lower bound becomes stable or after a certain number of iterations. As such, the \SDDiP algorithm combines these suggestions, stopping when either the lower bound has settled or an iteration limit is reached, whichever occurs first. 

In the Sampling step, a total of $M$ scenarios are compiled in $\Omega_I$ as realizations of the uncertain parameters $\xi^s_m$ for each stage in sequence. Next, the Forward step solves the appropriate forward problem from formulation~\eqref{SMCRP} for each stage $s \in \mathcal{S}\setminus\{0\}$ relative to the sampled scenario $m$ and characterized by the previous stage solution $x^{s-1, k}_{ij;Im}$. The Forward step obtains an updated feasible solution $\left(x^{s-1,k}_{ij;Im}, x^{sk}_{ij;m}, y^s_{tp;m} \right)$ as well as the objective function value across all decision-making stages $u_m$ for the corresponding scenario. Once the Forward step successfully updates the stage solution at all decision-making stages and collects the objective function values, the Upper bound update uses a Monte Carlo confidence interval approach to produce a statistically valid upper bound on the optimal objective value. Specifically, the Upper bound update uses the updated policy from the Forward step to compute the sample mean $\hat{\mu}$ and sample variance $\hat{\sigma}^2$ for use in computing a one-sided confidence value for the upper bound on line~14. The Upper bound update equation is a modification of the confidence interval given in Section~7 of Ref.~\cite{Fullner2025SDDPReview}, relying on the Central Limit Theorem $\frac{\hat{\mu} - \sum_{m=1}^M u_m}{\hat{\sigma}/\sqrt{M}} \approx \mathcal{N}(0, 1)$ for sufficiently large $M$ such that the probability of $\left(\sum_{m=1}^M u_m \le \hat{\mu} + z_{\alpha/2}\frac{\hat{\sigma}}{\sqrt{M}}\right) \approx 1 -\alpha$ applies the upper bound with confidence $1 - \alpha$. A common implementation of $z_{\alpha/2}$ is relative to a \SI{95}{\%} confidence upper bound, with a value of $z_{\alpha/2} = 1.96$ \cite{Pereira1991MSO}. Following the upper-bound update, the Backward step works through each stage in reverse order, intending to update the lower approximation utilizing cut coefficients $\left(v^s_{Imw}, \lambda^{sk}_{ij;Imw}\right)$ from suitable relaxations $R^s_{Im}$. Section 3.3 of Ref.~\cite{Zou2019SDDiP} details the sufficient conditions that a cut must satisfy for the \SDDiP algorithm to be valid, while this paper details two such cut families in the following paragraphs. Due to the nature of the Backward step and the cut families, the terminal stage $s=S$ will never accumulate cuts, so $\psi^s_{I}\left(\cdot\right) \equiv 0, \forall I \ge 0$. The Backward step continues until valid cuts have been aggregated for all decision-making stages $s \in \mathcal{S}\setminus\{0, S\}$ and proceeds to the final step of the algorithm. The final step of the \SDDiP algorithm is the lower-bound update, which solves the deterministic first-stage forward problem with updated lower approximation $\psi^1_{I+1}$, setting the optimal objective function value as the lower bound according to Eq.~\eqref{SMCRP:obj}. 

\begin{algorithm}[!ht] 
    \hypertarget{SDDiP}{}
    \DontPrintSemicolon
    \caption{\textcolor{myblue}{\textsf{SDDiP}}}
    \label{alg:SDDiP}
    Initialize: $LB \gets -\infty, ~ UB \gets +\infty, ~ I \gets 1, ~\{\psi^s_I(\cdot)\}_{s\in\mathcal{S}\setminus\{0\}} \gets \bm{0}$\;
    \While{$LB$ has not settled OR $I$ is less than an iteration limit}{
        /*Sampling step*/ \;
        Sample $M$ scenarios $\Omega_I = \{\xi^1_m, \xi^2_m, \ldots, \xi^S_m\}_{m = 1, \ldots, M}$\;
        /*Forward step*/\;
        \For{$m = 1, \ldots, M$}{
            \For{$s=1, \ldots, S$}{
                Solve forward problem $P^s_{I}\left(x^{s-1,k}_{ij;Im}, \psi^s_{I}, \xi^s_m\right)$ and collect solution $\left(x^{s-1,k}_{ij;Im}, x^{sk}_{ij;m}, y^s_{tp;m} \right)$\;
            }
            $u_m = \sum_{s\in\mathcal{S}\setminus\{0\}} \left(-\sum_{t\in\mathcal{T}^s} \sum_{p\in\mathcal{P}_m} \pi^s_{tpm} y^s_{tp;m} \right)$\;
        }
        /*(statistical) Upper bound update*/\;
        $\hat{\mu} \gets \frac{1}{M}\sum_{m=1}^M u_m$ and $\hat{\sigma}^2 \gets \frac{1}{M-1}\sum_{m=1}^M(u_m-\hat{\mu})^2$\;
        $UB \gets \hat{\mu} + z_{\alpha/2}\frac{\hat{\sigma}}{\sqrt M}$\;
        /*Backward step*/\;
        \For{$s=S,\ldots,2$}{
            \For{$m=1,\ldots,M$}{
                \For{$w=1,\ldots,W$}{
                    Solve a suitable relaxation $(R^s_{Im})$ of the updated problem $P^s_{I}\left(x^{s-1,k}_{ij;Im}, \psi^s_{I+1}, \xi^s_w\right)$ and collect cut coefficients $\left(v^s_{Imw}, \lambda^{sk}_{ij;Imw}\right)$\;
                }
                Add cuts using the coefficients $\left(v^s_{Imw}, \lambda^{sk}_{ij;Imw}\right)$ to $\psi^{s-1}_I$ to get $\psi^{s-1}_{I+1}$ 
            }
        }
        /*Lower bound update*/\;
        Solve $P^1_I\left(x^{0k}_{ij}, \psi^1_{I+1}\right)$ and set $LB$ to the optimal value\;
        $I \gets I+1$\;
    }
\end{algorithm}

Following the suggestions of the computational experiments of Ref.~\cite{Zou2019SDDiP}, two cut families that are proven to be applicable are used in the Backward step of the \SDDiP algorithm, those being the Integer Optimality Cuts and the Strengthened Benders' Cuts. Integer Optimality Cuts localize the cut information to a specific candidate solution and, as such, may not be as effective at other feasible solutions under nominal conditions if used in isolation. Alternatively, Strengthened Benders' Cuts rely less on localization to a candidate solution, thus affording additional computational benefits. As such, the use of these two cut families in cooperation is nearly as effective as other combinations explored in Ref.~\cite{Zou2019SDDiP}, with the benefit of a lower computational burden. 

The Integer Optimality Cuts are generated by evaluating the Forward step problems from formulation~\eqref{SMCRP} at feasible solutions. Specifically, the relaxation $R^{s+1}_{Iw}$ is the original \SMCRP problem $P^{s+1}_I\left(x^{sk}_{ij;Im}, \psi^{s+1}_{I+1}, \xi^{s+1}_w\right)$ at the candidate solution $x^{sk}_{ij;Im}$ obtained in the Forward step. As such, let $v^{s+1}_{I+1,mw}$ be the optimal objective function value of $P^{s+1}_I\left(x^{sk}_{ij;Im}, \psi^{s+1}_{I+1}, \xi^{s+1}_w\right)$ for a given $x^{sk}_{ij;Im}$ in all scenarios $w \in \mathcal{W}$. The Integer Optimality Cut added to $P^s_{I}\left(x^{s-1,k}_{ij;Im}, \psi^s_{I+1}, \xi^s_w\right)$ takes the following form:
\begin{equation}
    \theta^s \ge \left( \bar{v}^s_{I+1,m} - L^s \right) \left( \sum_{k\in\mathcal{K}} \sum_{i \in \mathcal{J}^{s-1,k}} \sum_{j \in \mathcal{J}^{sk}} \left( \bar{x}^{sk}_{ij;Im} - 1 \right)x^{sk}_{ij} + \sum_{k\in\mathcal{K}} \sum_{i \in \mathcal{J}^{s-1,k}} \sum_{j \in \mathcal{J}^{sk}} \left( x^{sk}_{ij} - 1 \right)\bar{x}^{sk}_{ij;Im} \right) + \bar{v}^s_{I+1,m}
\end{equation}
where $\bar{v}^s_{I+1,m} = \sum_{w\in\mathcal{W}}E_wv^{s+1}_{I+1,mw}$ is the expected objective function value localized to the solution $\bar{x}^{sk}_{ij;Im}$. In this case, the values of $(\bar{v}^s_{I+1,m}, \bar{x}^{sk}_{ij;Im})$ take the place of $(v^s_{Imw}, \lambda^{sk}_{ij;Imw})$ in line~17. Given the dependency of this cut family on the candidate solution, there is one cut added to the lower approximation for each sampled scenario $m = 1, \ldots, M$ for a total of $M$ Integer Optimality Cuts added to each stage in each iteration of the \SDDiP algorithm. 

The Strengthened Benders' Cuts are generated by evaluating the linear programming (LP) relaxation of the Forward step problems, as well as the Lagrangian relaxation of the Forward step problems in formulation~\eqref{SMCRP}. While the typical Benders' Cut relies only on the LP relaxation, wherein the cut coefficients $\left(v^s_{Imw}, \lambda^{sk}_{ij;Imw}\right)$ are the optimal objective function value of the LP relaxation and the basic optimal dual solution of constraints~\eqref{SMCRP:Link}, respectively, the Strengthened Benders' Cuts enhance the cut using the Lagrangian relaxation. The Strengthened Bender's Cut added to $P^s_{I}\left(x^{s-1,k}_{ij;Im}, \psi^s_{I+1}, \xi^s_w\right)$ takes the following form:
\begin{equation}
    \theta^s \ge \sum_{w\in\mathcal{W}} E_w\tilde{v}^{s+1}_{I+1,w} + \sum_{w\in\mathcal{W}}E_w\left( \sum_{k\in\mathcal{K}} \sum_{i \in \mathcal{J}^{s-1,k}} \sum_{j \in \mathcal{J}^{sk}} \tilde{\lambda}^{s+1, k}_{ij;I+1,w}x^{sk}_{ij} \right)
\end{equation}
where $\tilde{v}^s_{I+1,w} = \mathcal{L}^s_{I+1,w}\left( \tilde{\lambda}^{sk}_{ij;I+1,w} \right)$ is the optimal objective function of the Lagrangian relaxation, $\tilde{\lambda}^{sk}_{ij;I+1,w}$ are the basic optimal LP relaxation dual solution multipliers corresponding to constraints~\eqref{SMCRP:Link}, the LP relaxation is the \SMCRP with constraints~\eqref{SMCRP:x} relaxed to continuous on the range $[0, 1]$, and the Lagrangian relaxation is denoted as: 
\begin{subequations}
    \begin{alignat}{2}
\mathcal{L}^s_{Iw}\left(\tilde{\lambda}^{sk}_{ij;Iw}\right) \coloneq & \min \quad \left(-\sum_{t \in \mathcal{T}^s} \sum_{p \in \mathcal{P}_w} \pi^s_{tpw}y^s_{tp}\right) + \theta^s_I - \sum_{k\in\mathcal{K}} && \sum_{i \in \mathcal{J}^{s-1,k}} \sum_{j \in \mathcal{J}^{sk}} \tilde{\lambda}^{sk}_{ij;Iw}z^{sk}_{ij} 
\label{lag:obj} \\
& \text{s.t.} \quad \sum_{j \in \mathcal{J}^{sk}} x^{sk}_{ij} - \sum_{j^\prime \in \mathcal{J}^{s-2, k}} z^{sk}_{j^\prime i} = 0, \quad && \forall k \in \mathcal{K}, \forall i \in \mathcal{J}^{s-1, k} 
\label{lag:flow} \\
& \qquad \sum_{k \in \mathcal{K}} \sum_{i \in \mathcal{J}^{s-1,k}} \sum_{j \in \mathcal{J}^{sk}} V^{sk}_{tjpw} x^{sk}_{ij} \ge y^s_{tp}, \quad && \forall t \in \mathcal{T}^s, \forall p \in \mathcal{P}_w
\label{lag:vis}\\
& \qquad x^{sk}_{ij} \in \{0,1\}, \quad && \forall k \in \mathcal{K}, \forall i \in \mathcal{J}^{s-1, k}, \forall j \in \mathcal{J}^{sk} 
\label{lag:x} \\
& \qquad z^{sk}_{ij} \in [0,1], \quad && \forall k \in \mathcal{K}, \forall i \in \mathcal{J}^{s-2, k}, \forall j \in \mathcal{J}^{s-1,k} 
\label{lag:z} \\
& \qquad y^s_{tp} \in [0,1], \quad && \forall t \in \mathcal{T}^s, \forall p \in \mathcal{P}_w 
\label{lag:y} \\
& \qquad \theta^s \ge L^s 
\label{lag:LB}\\
& \qquad \theta^s \ge \text{cuts}_{I}
\label{lag:cuts}
    \end{alignat}
    \label{SMCRP:lag}
\end{subequations}
relative to each scenario $w \in \mathcal{W}$. Similarly to the Integer Optimality Cuts, one cut is added to the lower approximation for each sampled scenario as a result of the LP relaxation reliance on the previous stage solution $x^{s-1, k}_{ij;Im}$, for a total of $M$ Strengthened Benders' Cuts added to each stage in each iteration of the \SDDiP algorithm. Therefore, a total of $2M$ cuts are added in each iteration to improve the lower approximation and drive the \SDDiP algorithm to convergence. 

Reference~\cite{Zou2019SDDiP} details proofs of sufficient cut conditions in Section 3.3, as well as the proof of convergence to an optimal solution of an MSIP with binary state variables in finite iterations within Section 3.4. Therefore, since the \SMCRP is shown to match the MSIP format and obey the assumptions previously detailed, the \SDDiP converges to an optimal solution of the \SMCRP. Moreover, since the \SMCRP is an appropriate reformulation of the \MCRP, said solution returns an optimal policy of orbital maneuvers for each scenario.

\section{Computational Experiments} \label{sec:Experimentation}

To demonstrate the performance of the \SMCRP solved through the \SDDiP algorithm, two computational experiments are conducted using scenarios generated in unique ways. Several additional solution methods are explored for comparative purposes, including the use of non-maneuverable satellites, a random maneuver policy, and two dynamic programming methods to solve an MDP reformulation of the \SMCRP. The constellation of non-maneuverable satellites assumes that the satellites remain in their initial condition orbital slots for the duration of the time horizon in each distinct scenario. Separately, the random maneuver policy allows the constellation to select a random maneuver in each decision-making stage for each scenario and is evaluated \num{1000} times to determine the average performance. The solution of the \SMCRP through random orbital maneuvers is referred to as \textsf{Rand}. Finally, the deterministic upper bound is given by solving the \MCRP for each scenario in isolation. 

All experiments are conducted on a platform equipped with an Intel Core i9-12900 2.4 GHz (base frequency) CPU processor (16 cores and 24 logical processors) and 64 GB of RAM. Additionally, all experiments are conducted in Python (Version 3.9.6) \cite{Python}, with parameter generation of the VTWs $V^{sk}_{tjpw}$ conducted in MATLAB (Version R2024b, Update 6) with use of the Aerospace Toolbox \cite{MATLAB} and propagation using the Simplified General Perturbation 4 model. All optimization formulations, those being the \MCRP, the Forward \SMCRP problems, the LP relaxation of the \SMCRP problems, and the Lagrangian relaxation in formulation~\eqref{SMCRP:lag}, are solved using the commercial software package Gurobi Optimizer (Version 13.0.0) with default settings and a runtime limit of one hour. 

\subsection{Baseline Methods}

The \SMCRP applies to the Markovian nature of MDP formulations due to the removal of constraints~\eqref{MCRP:cost} and the consideration of stagewise independence, and as such, the \SMCRP may be reformulated into an MDP and subsequently solved. An MDP is comprised of the state space $\mathcal{N}$, action space $\mathcal{A}$, transition function $\mathcal{Q}(n^\prime \vert n, a) \in \mathbb{R}_+$, and reward function $\mathcal{R}(n, a) \in \mathbb{R}$, where $n, n^\prime \in \mathcal{N}$ are the current and next state, respectively, and $a \in \mathcal{A}$ is the action taken at state $n$. Additionally, the transition function refers to the probability that action $a$ taken in state $n$ will result in arrival at state $n^\prime$, and the reward function $\mathcal{R}(n, a)$ is the obtained reward given action $a$ is taken in state $n$. The objective of an MDP is to determine the best action to take in each state, referred to as the Policy, that maximizes the rewards of each state-action pair according to the reward function.  

The state space of the \SMCRP is the collection of all permutations of stages, constellation formations, and scenarios, where the constellation formations are the Cartesian powers of the satellite orbital slots. Specifically, the constellation formations are the orbital slot $j \in \mathcal{J}^{sk}$ chosen by satellite $k \in \mathcal{K}$, resulting in $K$ total choices of $J^{sk}$ total orbital slots, so the constellation formations take the form $\mathcal{X} = \left\{\left(j^1, \ldots, j^K\right) \colon j^k \in \left\{1, \ldots, J^{sk}\right\}\right\} = \left\{1, \ldots, J^{sk}\right\}^K$. As such, the state space is defined: 
\begin{equation}
    \mathcal{N} = \left\{ \left(0, \mathcal{J}^{0k}, 0\right), (s, x, w) \colon s \in \mathcal{S} \setminus \{0\}, x \in \mathcal{X}, w \in \mathcal{W} \right\}
\end{equation}
with total cardinality $|\mathcal{N}| = 1+S\left(J^K\right)W$. The additional state $n = \left(0, \mathcal{J}^{0k}, 0\right)$ represents the initial conditions of the constellation before the start of the time horizon. The state space represents the decision-making stage of the \SMCRP, the currently inhabited orbital slots for each satellite of the constellation, and the target scenario. 

Following logically from the state space, the action space of the \SMCRP is the selection of the constellation formation. In this way, actions $a = 1, \ldots, J^K$ correspond to the constellation formation in $\mathcal{X}$ such that the action space is defined: 
\begin{equation}
    \mathcal{A} = \mathcal{X} = \left\{\left(j^1, \ldots, j^K\right) \colon j^k \in \left\{1, \ldots, J^{sk}\right\}\right\}
\end{equation}
with total cardinality $|\mathcal{A}| = J^K$. The action space represents the decision made by the satellites in the constellation, mapping directly to the roll of $x^{sk}_{ij}$ in the \SMCRP. 

The reward function of the \SMCRP is the observation rewards obtained by the occupied constellation relative to the current scenario and stage of the state $n$. Given that state $n$ has the form $n = (s, x, w) \in \mathcal{N} \setminus \left\{\left(0, \mathcal{J}^{0k}, 0\right)\right\}$, state $n$ references stage $s \in \mathcal{S} \setminus \{0\}$, constellation formation $x \in \mathcal{X}$, and scenario $w \in \mathcal{W}$. The constellation indices $x =\left(j^1, \ldots, j^K\right)\in \mathcal{X}$ may be reformatted to a more familiar form $x^k_j \in \{0, 1\}, ~\forall k \in \mathcal{K}, ~\forall j \in \mathcal{J}^{sk}$ via the following: 
\begin{equation*}
    x^k_j = \begin{cases}
        1, & \text{if } j = j^k \\
        0, & \text{otherwise}
    \end{cases}, \quad \forall k \in \mathcal{K}, \forall j \in \mathcal{J}^{sk}
\end{equation*}
With $x^k_j$ obtained from state $n$, constraints~\eqref{SMCRP:vis} may be adapted to indicate whether visibility of a target is acquired via:
\begin{equation*}
    y^s_{tp} = \begin{cases}
        1, & \text{if } \sum_{k \in \mathcal{K}}\sum_{j \in \mathcal{J}^{sk}} V^{sk}_{tjpw}x^k_j \ge 1 \\
        0, & \text{otherwise}
    \end{cases}, \quad \forall t \in \mathcal{T}^s, \forall p \in \mathcal{P}_w
\end{equation*}
Then, utilizing the obtained visibility indicated by $y^s_{tp}$, the reward function is defined: 
\begin{equation}
    \mathcal{R}(n, a) = E_w\sum_{t \in \mathcal{T}^s} \sum_{p \in \mathcal{P}_w} \pi^s_{tpw}y^s_{tp}
\end{equation}
The reward function additionally incorporates the probability that scenario $w$ occurs through the use of $E_w$ and the target rewards $\pi^s_{tpw}$. Finally, the reward of the initial condition state is zero: $\mathcal{R}\left(\left(0, \mathcal{J}^{0k}, 0\right), a\right) = 0, ~\forall a \in \mathcal{A}$. 

The transition function of the \SMCRP is defined using the scenario probability of the next state $n^\prime$, additionally accounting for the action corresponding to the occupied constellation. Let the current state be $n = (s, x, w) \in \mathcal{N}$, the next state be $n^\prime = (s^\prime, x^\prime, w^\prime) \in \mathcal{N}$, and the action be $a \in \mathcal{A}$. The transition function is defined:
\begin{equation}
    \mathcal{Q}(n^\prime \vert n, a) = \begin{cases}
        E_{w^\prime}, & \text{if } s^\prime - s = 1 \text{ and } x^\prime = a \\
        1, & \text{if } s = S \text{ and } n^\prime = n \\
        0, & \text{otherwise}
    \end{cases}
\end{equation}
The consideration that transition function entries may equal one indicates that those states are terminal states, as there is no action to take in the final decision-making stage of reconfiguration when $s = S$. 

The MDP reformulation of the \SMCRP is solved in the computational experiments through the use of two well-known approaches, one of which is model-based, while the other is model-free. The use of these two approaches allows comparison when knowledge is and is not present while solving the MDP. The first solution method is value iteration, a dynamic programming offline planning algorithm \cite{Puterman1994MDPVI}. Value iteration assigns a value to each state-action pair $(n, a)$ and iteratively updates the value of each pair through the use of the Bellman operator \cite{Bellman2010DP} with the transition function until the variation between the value in the previous and current iterations is sufficiently low. Upon convergence, the most valuable action in each state is returned as a policy to maximize the rewards obtained within the reward function. The solution of the \SMCRP via value iteration is referred to as \textsf{VI}. The second solution method is Q-Learning, a model-free, value-based reinforcement learning algorithm \cite{Sutton1998Reinforcement}. Q-Learning also assigns a value to each state-action pair with iterative updates; however, Q-Learning balances exploration of random actions with exploitation of the current most optimal policy to determine the outcomes of various actions. As such, Q-Learning operates without direct knowledge of the transition function, learning the best action to take over iterations through intelligent trial and error. The Q-Learning process continues iterations until an iteration limit is reached and also returns a policy of which action to perform in each state. The solution of the \SMCRP via Q-Learning is referred to as \textsf{QL}. For all experiments, \textsf{VI} uses a discount of $0.99$ and a convergence tolerance of $0.01$; \textsf{QL} employs the same discount, a learning rate of $0.3$, a random action probability threshold of $0.99$, a decay factor of $0.999975$, and \num{1e6} Monte Carlo iterations.

\subsection{Experiment: Monitoring Orbital Targets}

An initial computational experiment is conducted concerning targets located in LEO as distinct scenarios alongside a constellation of satellites performing orbital maneuvers to gain observations. This experiment demonstrates the effectiveness of the \SDDiP algorithm concerning a large amount of target position uncertainty, as well as the effectiveness of orbital maneuverability to maintain visibility of highly mobile targets traveling at orbital velocities. Suppose that observation satellites are tasked to observe orbital targets in LEO that arrive at unique times within a mission horizon. Each target scenario is unique to ensure that no two scenarios have identical uncertain parameters, providing a variety of observation satellite orbits, orbital target orbits, and observation rewards, thus leading to a unique optimal orbital maneuver sequence in each scenario. 

\subsubsection{Design} 

\begin{figure}[!ht]
    \centering
    \includegraphics[width=0.5\textwidth]{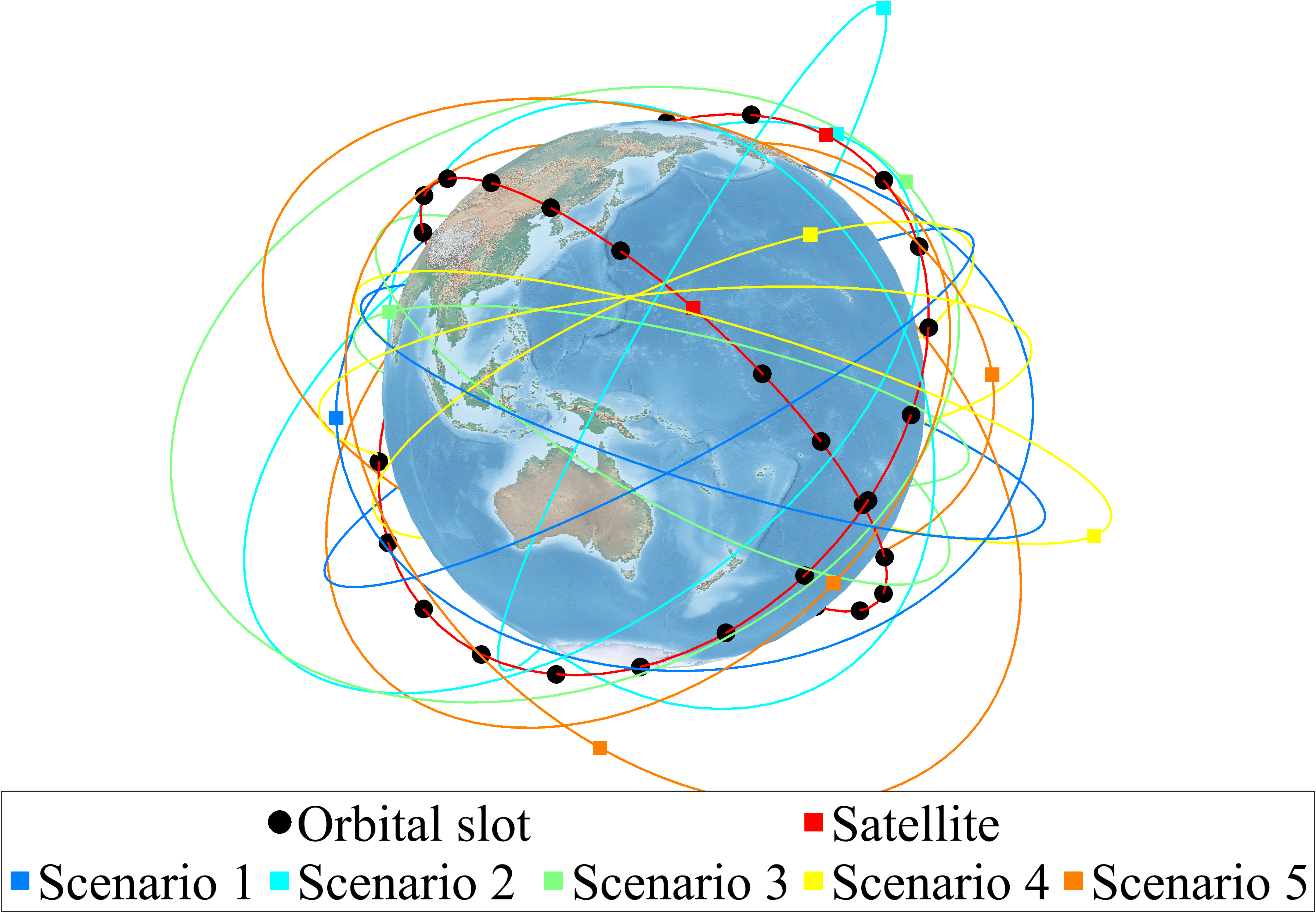}
    \caption{Orbits for monitoring orbital targets.}
    \label{fig:orbital_targets}
\end{figure}

The number of scenarios is set as $W = 20$ with $P_w = 10$ orbital targets in each scenario, and operations are conducted relative to $K = 2$ observation satellites. All $K$ observation satellites are assigned inclined circular orbits with an altitude between $500$ and $\SI{1000}{\kilo\meter}$, inclination between $40$ and $\SI{80}{\deg}$, and right ascension of the ascending node (RAAN) and argument of latitude between $0$ and $\SI{360}{\deg}$. Separately, all $P_w$ orbital targets are assigned inclined elliptical orbits with a periapsis altitude in the same range as the observation satellite's altitude, eccentricity between $0$ and $0.25$, inclination between $10$ and $\SI{80}{\deg}$, and RAAN, argument of periapsis, and true anomaly between $0$ and $\SI{360}{\deg}$. Additionally, $J^{sk} = 20$ orbital slots varying in argument of latitude are made available to each satellite $k \in \mathcal{K}$ in each stage $s\in\mathcal{S}\setminus\{0\}$, equally spaced between $0$ and $\SI{360}{\deg}$ with the inclusion of the initial phase. The orbits of the observation satellites, orbital slots, and a subset of five scenarios with three orbital targets from each are shown in Fig.~\ref{fig:orbital_targets}. The concept of the distribution of orbital targets depicted in Fig.~\ref{fig:orbital_targets} can be extrapolated to the full $W = 20$ scenarios with $P_w = 10$ orbital targets each, but is not shown directly due to the crowding that would result in the figure. 

The visibility of an orbital target $p\in\mathcal{P}_w$ from an observation satellite $k\in\mathcal{K}$ in orbital slot $j\in\mathcal{J}^{sk}$ is denoted by $V^{sk}_{tjpw}$ through the use of the intersatellite link condition presented in Ref.~\cite{Rogers2026ISL} and a maximum range condition. Let $\bm{r}^{sk}_{tj}$ and $\bm{r}^s_{tpw}$ be the position vectors of observation satellite $k$ in orbital slot $j$ and orbital target $p$ in scenario $w\in\mathcal{W}$, respectively, at time $t\in\mathcal{T}^s$ of stage $s\in\mathcal{S}\setminus\{0\}$. Additionally, let $\bm{\rho}^{sk}_{tpw}$ denote the relative position vector from said observation satellite to said orbital target at the same time and stage. Visibility of the orbital target from the observation satellite may be obtained only if the satellite is not occluded by Earth or Earth's atmosphere and is within range. From the geometry depicted in Ref.~\cite{Rogers2026ISL} with the addition of the ranging condition, the VTW of the orbital targets is denoted: 
\begin{equation*}
    V^{sk}_{tjpw} = \begin{cases}
        1, & \text{if } \sqrt{\left(r^{sk}_{tj}\right)^2 - (R_{\oplus} + \epsilon)^2} + \sqrt{\left(r^s_{tpw}\right)^2 - (R_{\oplus} + \epsilon)^2} - \rho^{sk}_{tjpw} > 0 \\
        & \text{and } \rho^{sk}_{tjpw} < \rho_{\max} \\
        0, & \text{otherwise}
    \end{cases}
\end{equation*}
where $R_{\oplus} = \SI{6378.14}{\kilo\meter}$ is the radius of the Earth, $\epsilon = \SI{100}{\kilo\meter}$ is an additional exclusion to compensate for atmospheric effects \cite{Bhattacherjee2019Altitude}, and $\rho_{\max} = \SI{1000}{\kilo\meter}$ is the maximum range for optical object detection \cite{Hussain2025Range}. 

Remaining parameters for the computational experiments are detailed as follows. The mission horizon is selected to be $4$ days beginning on June 1, 2026, at midnight (00:00:00) in Coordinated Universal Time (UTC). The mission horizon is discretized into $S=4$ equally spaced decision-making stages of orbital maneuvers. Furthermore, the discrete time step size $\Delta t$ is $100$ seconds, such that $T = 3456$ and $T^s = 864$. Separately, the observational rewards $\pi^s_{tpw}$ of each target are set to zero until a random time in a random stage, at which point they are set to one to model the orbital target ``appearing'' as a target for observation for the remaining duration of the mission. Finally, the initial lower bound of the lower approximation $\psi^s_{I}\left(\cdot\right)$ is set as $L^s = 0, \forall s \in \mathcal{S}\setminus\{0\}$. 

Due to the assumption of the \SDDiP algorithm that the initial conditions and first stage problem are both deterministic, a slight adaptation to the computational experiments is necessary to ensure a deterministic first stage problem. Specifically, an additional stage of decision-making is added before the start date of June 1, 2026, wherein the observation satellites are restricted from performing maneuvers or obtaining visibility of orbital targets. In this way, the mission operations begin at the start of the second decision-making stage (\textit{i.e.} the first reconfiguration stage) while the observation satellites maintain their initial conditions in the first decision-making stage with a deterministic objective function value of $0 + \psi^1_I\left(\cdot\right)$, thus incorporating the necessary cut information to update the lower bound and converge over iterations. 

\subsubsection{Results}

The results are detailed through the policy of each solution method, that being the sequence of orbital maneuvers performed by the constellation in each stage of reconfiguration $s\in\mathcal{S}\setminus\{0\}$. The policy of each solution method results in the return of the obtained observation rewards; the sum of the returned rewards for the entire mission is denoted $Z_w=\sum_{s\in\mathcal{S}\setminus\{0\}}\sum_{t\in\mathcal{T}^s}\sum_{p\in\mathcal{P}_w} \pi^s_{tpw}y^{s\ast}_{tpw},~\forall w\in\mathcal{W}$, where $y^{s\ast}_{tpw}$ is returned from:
\begin{equation*}
    y^{s\ast}_{tpw} = \begin{cases}
        1, & \text{if } \sum_{k \in \mathcal{K}}\sum_{i \in \mathcal{J}^{s-1,k}}\sum_{j \in \mathcal{J}^{sk}} V^{sk}_{tjpw} x^{sk\ast}_{ij} \ge 1 \\
        0, & \text{otherwise}
    \end{cases}, \quad \forall w \in \mathcal{W}, \forall s \in \mathcal{S} \setminus\{0\}, \forall t \in \mathcal{T}^s, \forall p \in \mathcal{P}_w
\end{equation*}
where $x^{sk\ast}_{ij}$ is the optimal policy of a given solution method. The following results are presented for analysis: 1) the sum of returned rewards $Z_w$ for all solution methods broken down by scenario and reconfiguration stage, 2) the improvement of the \SDDiP algorithm over all other solution methods, and 3) a comparison between the policy of the \SDDiP algorithm against the policy of the \MCRP deterministic upper bound for the scenario with the highest rewards. 

\begin{figure}[!ht]
    \centering
    \includegraphics[width=\textwidth]{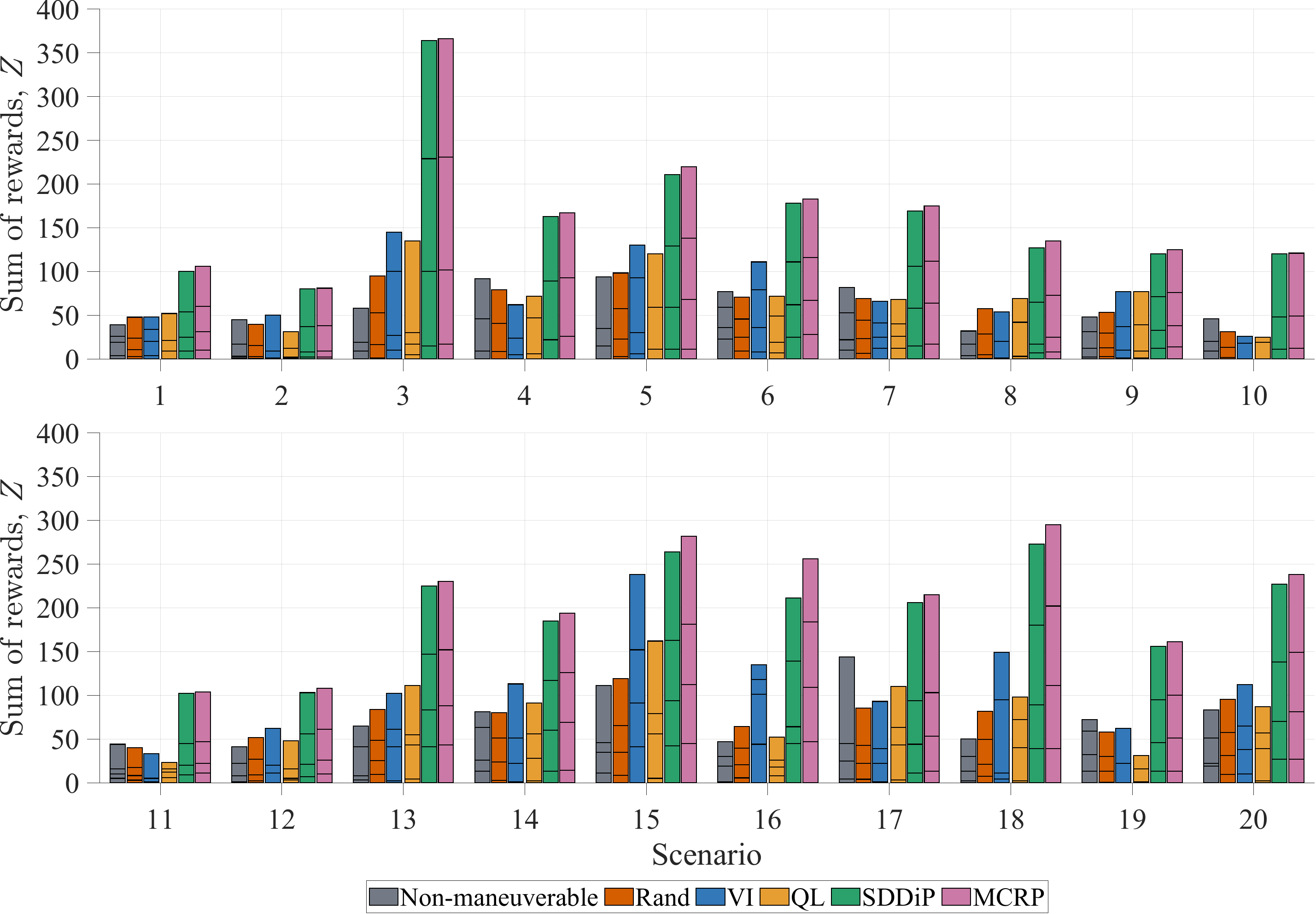}
    \caption{Sum of rewards $Z$ by scenario and stage for orbital monitoring.}
    \label{fig:Orbital_StageRewards}
\end{figure}

Figure~\ref{fig:Orbital_StageRewards} shows the sum of returned rewards, $Z_w$, as a bar chart for all solution methods in each scenario $w \in \mathcal{W}$. The figure additionally demonstrates the rewards obtained in each reconfiguration stage of the scenario, denoted by the horizontal lines within each bar, where stages $s=1, \ldots, S$ are ordered from bottom to top. Note that \textsf{Rand} reported values are the average performance of all $1000$ random maneuver evaluations. A clear trend is evident in the performance of each solution method: the \SDDiP algorithm is the highest-performing \SMCRP solution method in each scenario, while the non-maneuverable solution is the lowest-performing in $11$ scenarios, alongside \textsf{QL} in $4$, \textsf{Rand} in $3$, and \textsf{VI} in $2$. It can also be seen that the \SDDiP algorithm performs well, nearing the deterministic \MCRP upper bound in every scenario. While the \SDDiP algorithm guarantees global optimal convergence in finite iterations \cite{Zou2019SDDiP}, the algorithm was halted after $10$ iterations as a result of forward step problems exceeding the runtime limit of one hour set in the Gurobi optimizer. Despite halting before full convergence, the \SDDiP algorithm results shown in Fig.~\ref{fig:Orbital_StageRewards} outperform each other \SMCRP solution method by a wide margin in a majority of scenarios. Furthermore, while there are $11, ~ 5,$ and $9$ scenarios where the average \textsf{Rand} outperforms the non-maneuverable, \textsf{VI}, and \textsf{QL} solutions, respectively, the other scenarios perform much worse while said random maneuvers consume limited propellant resources. In addition, Table~\ref{tab:Orbital_Stats} depicts the minimum, maximum, average, and standard deviation of the sum of reward values shown in Fig.~\ref{fig:Orbital_StageRewards}. From these statistics, it is shown that the \SDDiP algorithm has the highest of every statistic out of the \SMCRP solution methods, while \textsf{QL} has the lowest minimum and \textsf{Rand} has the lowest maximum, mean, and standard deviation. 

\begin{table}[!ht]
    \centering
    \caption{Sum of reward statistics for orbital monitoring.}
    \begin{tabular}{ l  r r r r r r }
\hline \hline 
Statistic          & Non-maneuverable & \multicolumn{4}{c}{\SMCRP methods}           & \MCRP \\ 
\cmidrule(lr){3-6}
                   &                  & \textsf{Rand} & \textsf{VI} & \textsf{QL} & \SDDiP &       \\
\hline
Minimum            & 32.00  & 31.16  & 26.00  & 23.00  & 80.00  & 81.00  \\
Maximum            & 144.00 & 199.02 & 238.00 & 162.00 & 364.00 & 366.00 \\
Mean               & 67.55  & 70.02  & 93.40  & 76.70  & 179.20 & 188.10 \\
Standard deviation & 28.37  & 22.90  & 50.53  & 38.08  & 71.01  & 75.01  \\
\hline \hline 
    \end{tabular}
    \label{tab:Orbital_Stats}
\end{table}

\begin{figure}[!ht]
    \centering
    \includegraphics[width=\textwidth]{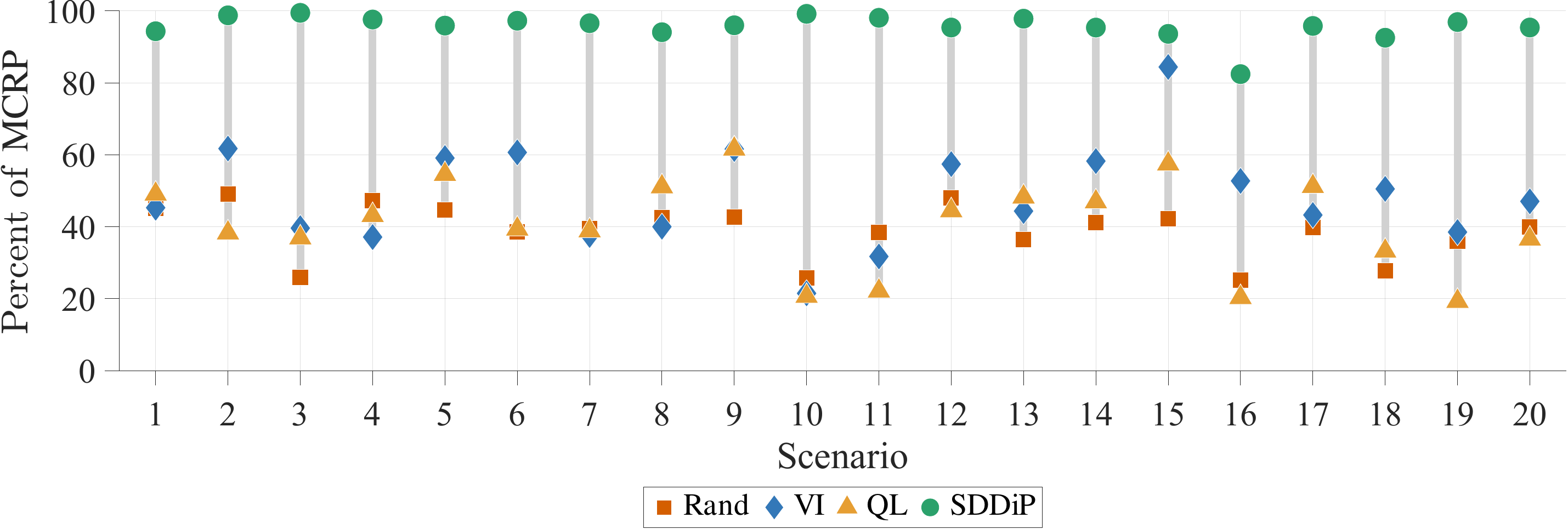}
    \caption{Improvement of \textcolor{myblue}{\textsf{SDDiP}} over other \textcolor{myblue}{\textsf{SMCRP}} solution methods for orbital monitoring.}
    \label{fig:Orbital_Improv_SDDiP}
\end{figure}

To further indicate the level of performance of each \SMCRP solution method, Fig.~\ref{fig:Orbital_Improv_SDDiP} depicts the relative sum of rewards as a percentage of the \MCRP deterministic upper bound. That is, each entry denotes the performance as a fraction of the \MCRP upper bound, denoting how close each solution method is to reaching the deterministic upper bound performance. As such, Fig.~\ref{fig:Orbital_Improv_SDDiP} shows that the \SDDiP algorithm outperforms every other \SMCRP solution method, being within \SI{5}{\%} of the \MCRP solution in $15$ scenarios. Table~\ref{tab:Orbital_Improv_Stats} additionally reports the statistics of the percent improvement of the \SDDiP algorithm over other methods, reporting that the \SDDiP algorithm outperforms every other \SMCRP solution method by more than \SI{100}{\%} on average as a result of the optimization-based approach. The \SDDiP algorithm has the lowest minimum and average improvement relative to \textsf{VI}, and the lowest maximum and standard deviation improvement relative to \textsf{Rand}, as well as the highest minimum relative to \textsf{Rand} and the highest maximum, average, and standard deviation relative to \textsf{QL}. As such, \textsf{VI} is the solution method with the solution quality closest to the \SDDiP algorithm, wherein Table~\ref{tab:Orbital_Stats} reflects the second-highest average sum of rewards and Table~\ref{tab:Orbital_Improv_Stats} reflects the lowest average improvement of the \SDDiP algorithm over \textsf{VI}. Logically, \textsf{VI} would perform more adequately than \textsf{QL} as a result of the knowledge of the transition function contained within the Bellman operator. 

\begin{table}[!ht]
    \centering
    \caption{Statistics of \textcolor{myblue}{\textsf{SDDiP}} solution relative to other \textcolor{myblue}{\textsf{SMCRP}} solution methods for orbital monitoring.}
    \begin{tabular}{ l  r r r r }
\hline \hline 
Statistic          & \multicolumn{3}{c}{\SMCRP methods}  & \MCRP \\ 
\cmidrule(lr){2-4}
                   & \textsf{Rand} & \textsf{VI} & \textsf{QL} &       \\
\hline
Minimum            & 98.70  & 10.92  & 55.84  & -17.58 \\
Maximum            & 285.16 & 361.54 & 403.23 & -0.55  \\
Mean               & 156.51 & 114.97 & 165.03 & -4.38  \\
Standard deviation & 57.10  & 75.92  & 106.47 & 3.62   \\
\hline \hline 
    \end{tabular}
    \label{tab:Orbital_Improv_Stats}
\end{table}

For a more thorough comparison between the \SDDiP algorithm solution and the \MCRP deterministic upper bound, Scenario $3$ is explored as an in-depth analysis, a selection made as a result of the large sum of rewards and gap between the \SMCRP solution methods. Table~\ref{tab:Orbital_varx} reports the orbital maneuver sequence of each satellite at all reconfiguration stages, as well as the rewards obtained by the satellite during that stage, including the orbital slot index ($j$) and the argument of latitude (Arg. Lat.). The overall sequence of orbital maneuvers is similar, with some key differences of a single orbital slot index in a single stage of each satellite maneuver sequence. Specifically, in Stage $1$, Satellite $2$ maneuvers into orbital slot $7$ in the \MCRP sequence and orbital slot $8$ in the \SDDiP sequence, resulting in two fewer rewards obtained due to a shift in target visibility from the resultant Arg. Lat. shift of \SI{18}{\deg}. Additionally, in Stage $2$, Satellite $1$ maneuvers into orbital slot $14$ in the \MCRP sequence and orbital slot $13$ in the \SDDiP sequence; however, the obtained rewards do not shift in this case. Upon closer inspection, the accumulation of rewards shifts the distribution, wherein the \MCRP maneuver obtains all $26$ rewards from a single target, and the \SDDiP maneuver obtains the rewards from three distinct targets. As reflected in Table~\ref{tab:Orbital_varx}, there are no maneuver sequence differences in Stages $3$ and $4$, resulting in identical rewards throughout the rest of this scenario. Additionally, more rewards are obtained in later stages of the scenario, a result given the assignment of rewards such that they are zero until a random time in a random stage, at which point the target ``appears'' for observation. As such, only one target is active in Stage $1$ approximately \SI{75}{\%} through the stage, which is joined by five more targets more than halfway through Stage $2$, while two more targets appear during Stage $3$, and the final two targets join at the start of Stage $4$. 

\begin{table}[!ht]
    \centering
    \caption{Orbital maneuver sequence of the \textcolor{myblue}{\textsf{MCRP}} and \textcolor{myblue}{\textsf{SDDiP}} for orbital monitoring of Scenario $3$.}
    \begin{tabular}{ l l  r r r  r r r}
\hline \hline
~ & ~ & \multicolumn{3}{c}{\MCRP} & \multicolumn{3}{c}{\SDDiP}\\
\cmidrule(lr){3-5} \cmidrule(lr){6-8}
Satellite, $k$ & Stage, $s$ & Index, $j$ & Arg. Lat., $\deg$ & Rewards & Index, $j$ & Arg. Lat., $\deg$ & Rewards \\
\hline
1 & 1 & 12 & 265.11 & 10.00 & 12 & 265.11 & 10.00 \\
  & 2 & 14 & 257.62 & 26.00 & 13 & 239.62 & 26.00 \\
  & 3 & 16 & 250.13 & 39.00 & 16 & 250.13 & 39.00 \\
  & 4 & 18 & 242.64 & 41.00 & 18 & 242.64 & 41.00 \\
2 & 1 & 7  & 265.98 & 7.00  & 8  & 283.98 & 5.00  \\
  & 2 & 10 & 161.75 & 59.00 & 10 & 161.75 & 59.00 \\
  & 3 & 11 & 21.53  & 90.00 & 11 & 21.53  & 90.00 \\
  & 4 & 13 & 259.30 & 94.00 & 13 & 259.30 & 94.00 \\
\hline \hline
    \end{tabular}
    \label{tab:Orbital_varx}
\end{table}

Overall, the case of monitoring orbital targets via LEO observation satellites highlights the performance of the \SDDiP algorithm in a highly varied set of scenarios. The \SDDiP algorithm results approach the \MCRP deterministic upper bound in most scenarios, while outperforming the other \SMCRP solution methods in every scenario by a wide margin on average. The high level of performance originates from the optimal policy of the \SDDiP algorithm closely mirroring the orbital maneuver sequence of the \MCRP, often differing by singular orbital slots sparsely throughout the mission duration. Additionally, while not fully converged to the global optimal policy due to an imposed runtime limit on the optimization software, the \SDDiP algorithm performs consistently while optimizing the policy for each scenario simultaneously within Algorithm~\ref{alg:SDDiP}. An additional computational experiment is conducted under different scenario generation techniques to further demonstrate the versatility of the \SDDiP algorithm. 

Although the \SDDiP algorithm outperforms all other \SMCRP solution methods in each scenario, there are important tradeoffs that should be noted. Firstly, the computation time of each solution method is crucial. The MDP solution methods, as heuristics, achieve convergence or reach an iteration limit within less than $12$ hours. The \SDDiP algorithm alternatively relies on the \SMCRP dynamic programming equations in formulation~\eqref{SMCRP}, as well as the relaxations for the collection of cut coefficients, which are mainly MILP formulations. As such, the $10$ iterations of the \SDDiP algorithm complete after a total of 60 hours, more than five times the convergence time of the slowest MDP solution method. Secondly, the MDP solution methods are highly sensitive to their hyperparameters, especially the discount factor, convergence threshold, and Monte Carlo iteration limit. This work does not consider a hyperparameter search to tune these values and increase the effectiveness of the MDP solution methods, which may increase their performance. Therefore, given these tradeoffs, the \SDDiP outperforms all other \SMCRP solution methods under the given experiment parameters.

\subsection{Experiment: Monitoring Hurricanes}

In addition to the computational experiment concerning orbital targets, simulated hurricane paths are utilized as target sets to demonstrate the effectiveness of the \SDDiP algorithm relative to a realistic use case concerning the highly dynamic nature of hurricane natural disasters. The observation satellites are referred to as satellites in the context of this experiment for brevity. Each simulated hurricane is a distinct scenario and is simulated with randomness to ensure that no two hurricane paths are identical, providing a variety in overall path trajectories. 

\subsubsection{Design}

The number of simulated hurricanes, and thus scenarios, is similarly set as $W = 20$ and is made up of $P_w = 16$ discrete target points along the hurricane path. As in previous literature, the time between hurricane points is considered to be $6$ hours \cite{Pearl2025Benchmarking,Pearl2026REOSSP}, for a mission horizon of $4$ days. Many parameters utilized in this experiment remain the same as in the experiment concerning orbital targets, with some key exceptions. For instance, the satellite constellation's initial conditions remain unchanged from Fig.~\ref{fig:orbital_targets}, and the mission start time, discrete time step size, and initial lower bound of the lower approximation also remain unchanged. 

\begin{figure}[!ht]
    \centering
    \begin{subfigure}[h]{0.49\textwidth}
        \centering
        \includegraphics[width = 0.75\textwidth]{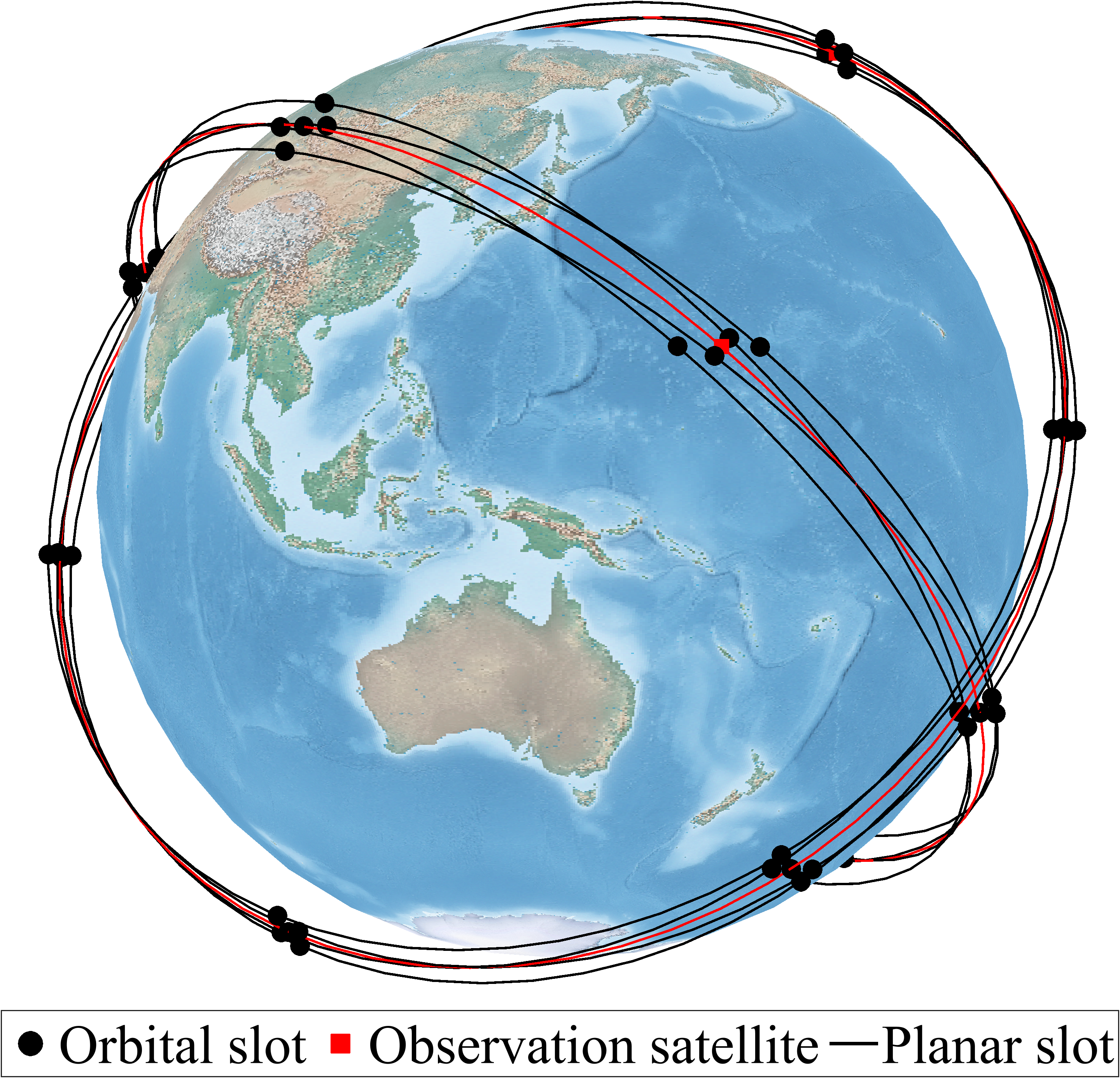}
        \caption{Orbital slots with planar options.}
        \label{fig:PlanarSlots}
    \end{subfigure}
    \hfill
    \begin{subfigure}[h]{0.49\textwidth}
        \centering
        \includegraphics[width = \textwidth]{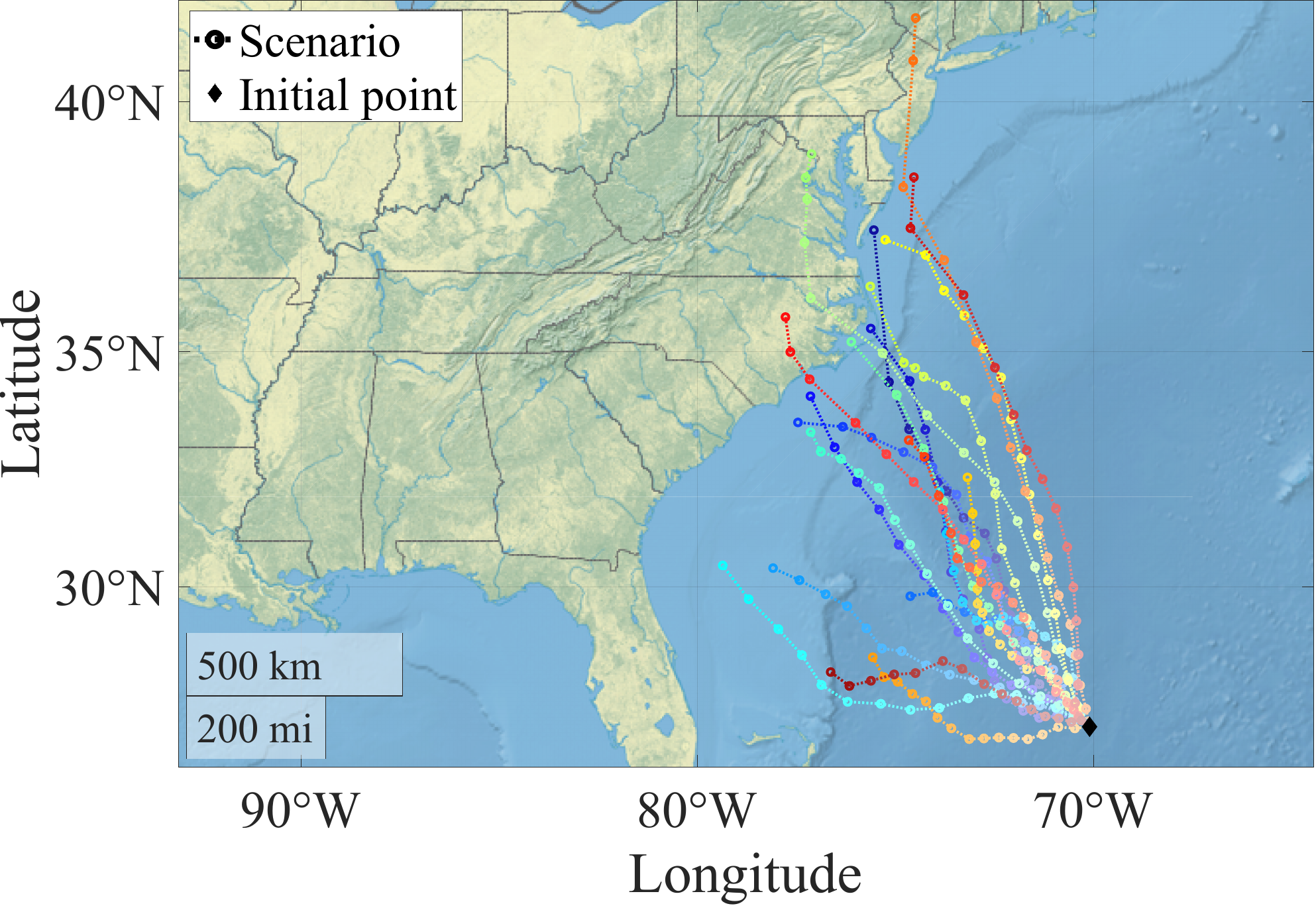}
        \caption{Hurricane scenarios.}
        \label{fig:TC_targets}
    \end{subfigure}
    \caption{Satellite and Hurricane parameters.}
    \label{fig:TC_Params}
\end{figure}

In this computational experiment, the orbital slots are modified to include changes in the plane of the orbital slot, as well as the argument of latitude, in the same manner as Ref.~\cite{Pearl2026REOSSP}. Specifically, the orbital slots include changes in inclination, RAAN, and argument of latitude, where plane changes occur on either side $(\pm)$ of the initial condition plane. In this experiment, a total of three planes of inclination and RAAN are utilized, resulting in a total of five unique plane combinations, with an additional six argument of latitude slots within each plane, for a total of $J^{sk} = 30$ unique orbital slots. The distribution of the orbital slots is depicted in Fig.~\ref{fig:PlanarSlots}. As a result of the new orbital slots, a total of seven different maneuvers are possible, those being changes in: 1) argument of latitude, 2) inclination, 3) RAAN, 4) inclination and RAAN, 5) inclination followed by argument of latitude, 6) RAAN followed by argument of latitude, and 7) inclination and RAAN followed by argument of latitude. Due to the slightly more expensive nature of plane change maneuvers, $S=2$ decision-making stages are considered in this experiment. Therefore, the number of time steps in the mission horizon and each stage is $T = 3456$ and $T^s = 1728$, respectively. 

The various scenarios of simulated hurricanes are propagated using a Markov-chain model from Ref.~\cite{Cui2019TCSim} with induced Brownian motion modified using an AR(1) autoregressive framework \cite{Iacovacci2016AR1}. The Markov-chain model utilized to propagate each hurricane path is defined as follows:
\begin{equation}
    \begin{bmatrix} g^{\text{lat}}_{p+1} \\ g^{\text{lon}}_{p+1} \end{bmatrix} = 
    \begin{bmatrix} g^{\text{lat}}_p \\ g^{\text{lon}}_p \end{bmatrix} + 
    \begin{bmatrix} \mu^{\text{lat}}\left( g^{\text{lat}}_p, g^{\text{lon}}_p \right) \\ \mu^{\text{lon}}\left( g^{\text{lat}}_p, g^{\text{lon}}_p \right) \end{bmatrix}\delta p + 
    \begin{bmatrix} \sigma^{\text{lat}}\left( g^{\text{lat}}_p, g^{\text{lon}}_p \right) & 
    \sigma^{\text{lat,lon}}\left( g^{\text{lat}}_p, g^{\text{lon}}_p \right) \\ 
    \sigma^{\text{lat,lon}}\left( g^{\text{lat}}_p, g^{\text{lon}}_p \right) & 
    \sigma^{\text{lon}}\left( g^{\text{lat}}_p, g^{\text{lon}}_p \right)\end{bmatrix}\sqrt{\delta p}Z_p, \quad \forall p \in \mathcal{P}_w\setminus\{P_w\}
    \label{eq:TC}
\end{equation}
where $\left(g^{\text{lat}}_p, g^{\text{lon}}_p \right)$ is the latitude longitude pair of target $p \in \mathcal{P}_w$ in the hurricane path and $\delta p$ is the $6$ hour step size between pairs. Additionally, $\mu^{\text{lat}}$ and $\mu^{\text{lon}}$ are the drift terms that represent the hurricane's acceleration component, and $\sigma^{\text{lat}}$ and $\sigma^{\text{lon}}$ are the diffusion terms that represent the randomness of the acceleration according to the random variable $Z_p$. The additional term $\sigma^{\text{lat,lon}}$ is the correlation between the two directions expressed as the root-covariance $r\left(g^{\text{lat}}_p, g^{\text{lon}}_p\right)\sqrt{\sigma^{\text{lat}}\left(g^{\text{lat}}_p, g^{\text{lon}}_p\right)\sigma^{\text{lon}}\left(g^{\text{lat}}_p, g^{\text{lon}}_p\right)}$. While Ref.~\cite{Cui2019TCSim} considers simple Brownian motion for the random variable $Z_p$ without memory, this paper adopts an AR(1) autoregressive framework to ensure that each hurricane scenario deviates further from a set of initial conditions as it progresses along the hurricane path. The progression of $Z_p$ is governed by: 
\begin{equation*}
    Z_{p+1} = \rho Z_p + \sqrt{1 - \rho^2}\mathcal{N}(0, 1), \quad \forall p \in \mathcal{P}_w\setminus\{P_w\}
\end{equation*}
utilizing a scaling factor $\rho = 0.85$ and a random variable from the normal distribution $\mathcal{N}(0, 1)$ \cite{Iacovacci2016AR1}. It is assumed that $Z_1 = 0$, while the values of the drift, diffusion, and correlation functions are determined from Figs.~4 and~5 of Ref.~\cite{Cui2019TCSim}. Using the initial conditions $\left(g^{\text{lat}}_1, g^{\text{lon}}_1 \right)$ matching the starting point of Hurricane Joaquin, the scenario hurricane paths are shown in Fig.~\ref{fig:TC_targets}. Hurricane Joaquin occurred in 2015 and heavily impacted the Bahamas and South Carolina, causing $34$ deaths and \$60 million in damage directly \cite{Berg2016Joaquin} as well as $19$ deaths and \$12 billion due to flooding \cite{Fritz2016SouthCarolina}, while reaching Category Four status on October 1, 2015, at 9:00 PM UTC with peak winds of $115$ mph. The starting point of Hurricane Joaquin is considered to be the first occurrence of tropical storm status (wind speed of at least $39$ mph \cite{TempestDefinitions_2022}) obtained from Ref.~\cite{WUnderground}. 

Since each target point $p \in \mathcal{P}_w$ of each hurricane scenario is a distinct point along the scenario path, the time-dependent rewards $\pi^s_{tpw}$ are used to model the motion of the hurricane over time. Firstly, visibility $V^{sk}_{tjpw} = 1$ if the target falls within a conical field of view of \SI{70}{\deg} of the relevant satellite at the relevant time step and stage, and $V^{sk}_{tjpw} = 0$ otherwise. Secondly, the time-dependent rewards represent the existence of the target $p$ at time $t\in\mathcal{T}^s$ of stage $s\in\mathcal{S}\setminus\{0\}$. As such, the rewards follow: 
\begin{equation*}
    \pi^s_{tpw} = \begin{cases}
        1, & \text{if } 1 + (p-1)(T/P_w) \le t \le p(T/P_w)\\
        0, & \text{otherwise}
    \end{cases}, \quad \forall s \in \mathcal{S}\setminus\{0\}, \forall t \in \mathcal{T}^s, \forall p \in \mathcal{P}_w, \forall w \in \mathcal{W}
\end{equation*}

\subsubsection{Results}

\begin{figure}[!ht]
    \centering
    \includegraphics[width=\textwidth]{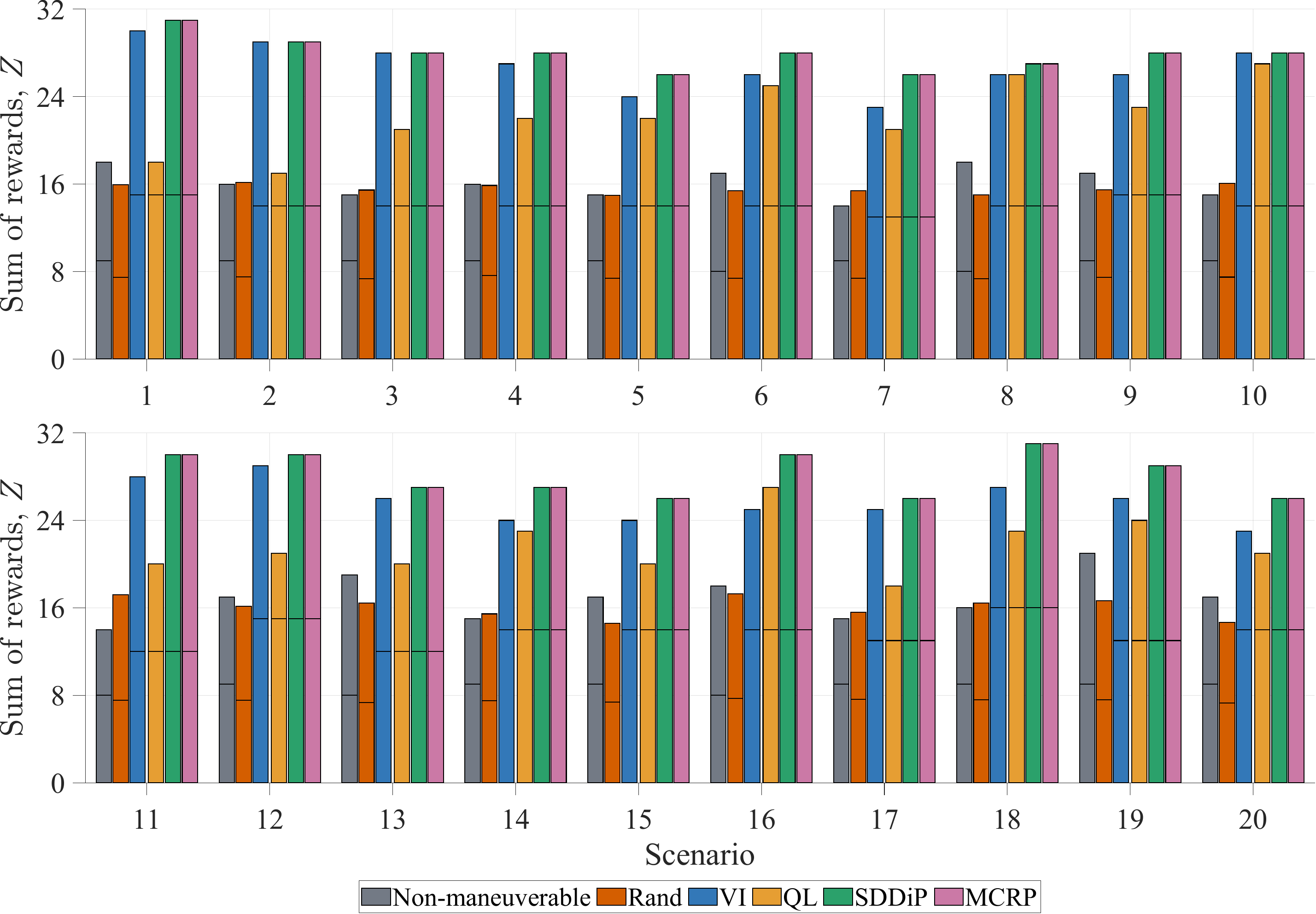}
    \caption{Sum of rewards $Z$ by scenario and stage for hurricane monitoring.}
    \label{fig:TC_StageRewards}
\end{figure}

The results of the computational experiment relative to dynamic hurricanes reflect similar conclusions to those concerning orbital target monitoring, with a slightly reduced scale of total rewards obtained due to the nature of the encoded rewards and available visibility. Figure~\ref{fig:TC_StageRewards} shows the sum of returned rewards, $Z_w$, as a bar chart, with a single horizontal bar bisecting the first stage from the second. As with the previous experiment, the \SDDiP algorithm is the highest-performing in each scenario by achieving the deterministic upper bound of the \MCRP, showing full convergence of the \SDDiP algorithm within a total of $10$ iterations (for an applicable comparison with the previous experiment). Otherwise, the non-maneuverable solution is the lowest performing in $8$ scenarios, alongside \textsf{Rand} in the other $12$ scenarios, indicating that \textsf{VI} and \textsf{QL} outperform the non-maneuverable and \textsf{Rand} in each scenario. Furthermore, \textsf{VI} outperforms \textsf{QL} in all but two cases, where the methods tie in Scenario $8$, and \textsf{QL} outperforms \textsf{VI} in Scenario $16$. Additionally, Table~\ref{tab:TC_Stats} depicts the statistics reflected in Fig.~\ref{fig:TC_StageRewards} similarly to Table~\ref{tab:Orbital_Stats}. Of course, the \SDDiP algorithm has the highest statistics tied with the deterministic upper bound, while among the \SMCRP solution methods, \textsf{Rand} has the lowest statistics. This further reflects the sub-optimality of performing random orbital maneuvers and the effectiveness of the \SDDiP algorithm in determining an optimal reconfiguration policy. 

\begin{table}[!ht]
    \centering
    \caption{Sum of reward statistics for hurricane monitoring.}
    \begin{tabular}{ l  r r r r r r }
\hline \hline 
Statistic          & Non-maneuverable & \multicolumn{4}{c}{\SMCRP methods}           & \MCRP \\ 
\cmidrule(lr){3-6}
                   &                  & \textsf{Rand} & \textsf{VI} & \textsf{QL} & \SDDiP &       \\
\hline
Minimum            & 14.00 & 14.58 & 23.00 & 17.00 & 26.00 & 26.00 \\
Maximum            & 21.00 & 17.28 & 30.00 & 27.00 & 31.00 & 31.00 \\
Mean               & 16.50 & 15.81 & 26.20 & 21.95 & 28.05 & 28.05 \\
Standard deviation & 1.76  & 0.75  & 2.04  & 2.86  & 1.70  & 1.70  \\
\hline \hline 
    \end{tabular}
    \label{tab:TC_Stats}
\end{table}

The main indication as to why the \SDDiP algorithm converges in this computational experiment and not in the orbital target computational experiment lies within the differing structure of the problems in the Forward and Backward steps. In the orbital target experiment, every target has active obtainable rewards by the conclusion of the mission horizon, leading to a large number of obtainable rewards and potentially high visibility due to overlapping orbital alignment. However, concerning the modeling of the hurricane motion, only one target has active obtainable rewards at any given time step within the mission horizon, and the satellites only obtain visibility of each target upon a flyover of the geographic region, both of which limit the overall obtained sum of rewards. Therefore, while many orbital maneuver decisions may perform similarly with such a surplus of rewards in the orbital monitoring experiment, the optimal maneuver sequence relative to hurricanes is more limited. This aspect leads to a more favorable set of problem geometries in the search space for each subproblem of the \SDDiP algorithm, thus leading to a more rapid convergence to the global optimal.

\begin{figure}[!ht]
    \centering
    \includegraphics[width=\textwidth]{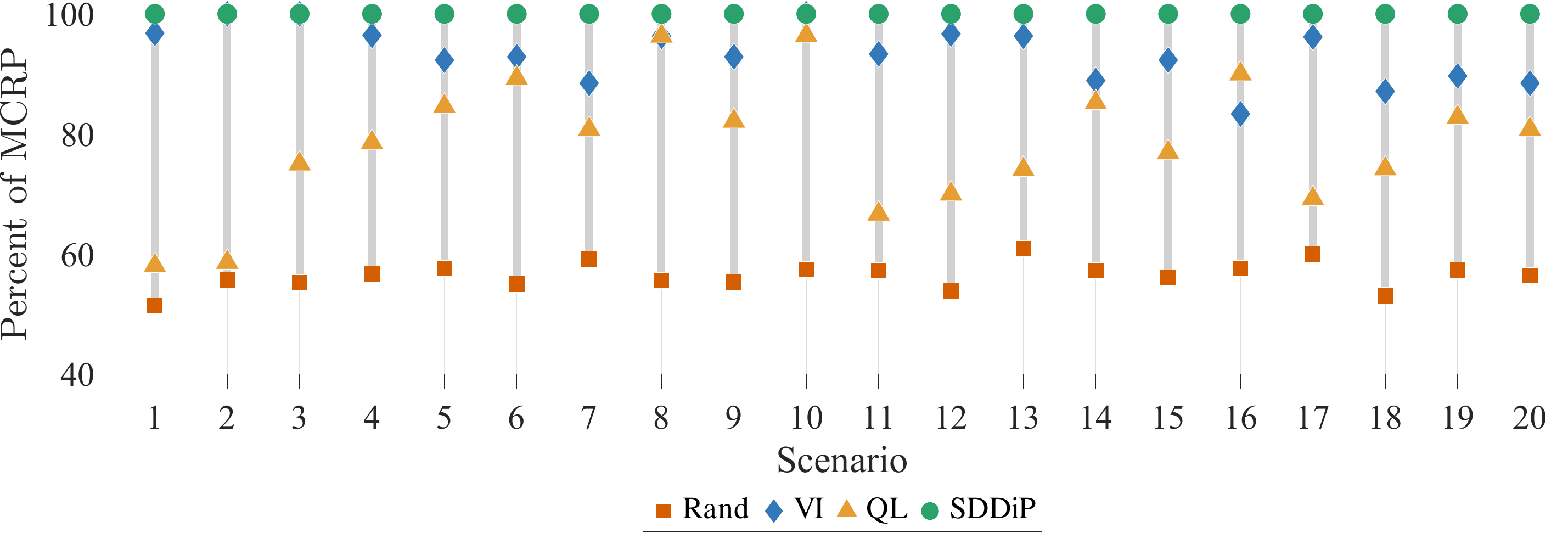}
    \caption{Improvement of \textcolor{myblue}{\textsf{SDDiP}} over other \textcolor{myblue}{\textsf{SMCRP}} solution methods for hurricane monitoring.}
    \label{fig:TC_Improv_SDDiP}
\end{figure}

Further exemplifying the convergence of the \SDDiP algorithm, Fig.~\ref{fig:TC_Improv_SDDiP} visualizes the percent of the \MCRP deterministic upper bound that each \SMCRP solution method achieves, reflecting that the \SDDiP algorithm achieves \SI{100}{\%} of the observation rewards also obtained by the \MCRP. Due to the lower total sum of rewards and resulting minimal gap in obtained rewards from the problem structure, the performance between each method is much closer than in the orbital target experiment. Specifically, the \SDDiP algorithm has a minimum improvement over \textsf{VI} of \SI{0}{\%} in Scenarios $2,~3,$ and $10$, meaning that the sum of observation rewards obtained by both solution methods is identical for these scenarios, additionally matching that of the \MCRP. Similarly, while the \SDDiP algorithm outperforms \textsf{Rand} the most in each statistic, the \SDDiP algorithm outperforms \textsf{VI} the least in each statistic, as reported in Table~\ref{tab:TC_Improv_Stats}. Again, \textsf{VI} outperforms \textsf{QL} due to the knowledge of the transition function, while the ability of policies obtained by \textsf{QL} has the capacity to outperform \textsf{VI} under the proper conditions, as indicated through the performance of \textsf{QL} in Scenario $16$. 

\begin{table}[!ht]
    \centering
    \caption{Statistics of \textcolor{myblue}{\textsf{SDDiP}} solution relative to other \textcolor{myblue}{\textsf{SMCRP}} solution methods for hurricane monitoring.}
    \begin{tabular}{ l  r r r r }
\hline \hline 
Statistic          & \multicolumn{3}{c}{\SMCRP methods}  & \MCRP \\ 
\cmidrule(lr){2-4}
                   & \textsf{Rand} & \textsf{VI} & \textsf{QL} &       \\
\hline
Minimum            & 64.23 & 0.00  & 3.70  & 0.00 \\
Maximum            & 94.53 & 20.00 & 72.22 & 0.00 \\
Mean               & 77.46 & 7.32  & 29.85 & 0.00 \\
Standard deviation & 7.12  & 5.47  & 18.99 & 0.00 \\ 
\hline \hline 
    \end{tabular}
    \label{tab:TC_Improv_Stats}
\end{table}

While the obtained rewards of the \SDDiP algorithm and \MCRP are identical in every scenario, Scenario $18$ is explored in-depth, due again to the large sum of rewards and gap between the \SMCRP solution methods. Table~\ref{tab:TC_varx} reports the orbital maneuver sequence of each satellite in the same manner as Table~\ref{tab:Orbital_varx}, additionally including the inclination (inc.) and RAAN of the orbital slot. Only a single maneuver differs within the optimal maneuver sequence, a proportional comparison due to there being only two stages of reconfiguration present. Specifically, in Stage $2$, Satellite $2$ maneuvers into orbital slot $19$ in the \MCRP sequence and orbital slot $12$ in the \SDDiP sequence, which changes the distribution of obtained rewards without changing the sum of obtained rewards. Specifically, due to the nature of the encoded hurricane motion within the rewards, target points $p = 9, \ldots, 16$ have obtainable rewards within Stage $2$, wherein each target has a total of $T/P = 216$ time steps in the provided order. As such, the \MCRP sequence obtains the rewards ordered 2--0--2--0--0--3--0--0, while the \SDDiP sequence obtains the rewards ordered 0--2--2--0--0--2--1--0, directly a result of the \SI{60}{\deg} shift in the argument of latitude between orbital slots $19$ and $12$. Notably, both satellites perform a change in inclination to be better positioned for observation in Stage $2$, while Satellite $2$ additionally performs a change in RAAN. 

\begin{table}[!ht]
    \centering
    \caption{Orbital maneuver sequence of the \textcolor{myblue}{\textsf{MCRP}} and \textcolor{myblue}{\textsf{SDDiP}} for hurricane monitoring of Scenario $18$.}
    \resizebox{\textwidth}{!}{
    \begin{tabular}{ l l  r r r r r  r r r r r }
\hline \hline
~ & ~ & \multicolumn{5}{c}{\MCRP} & \multicolumn{5}{c}{\SDDiP}\\
\cmidrule(lr){3-7} \cmidrule(lr){8-12}
Satellite, & Stage, & Index, & Inc,   & RAAN,  & Arg. Lat., & Rewards & Index, & Inc,   & RAAN,  & Arg. Lat., & Rewards \\
$k$        & $s$    & $j$    & $\deg$ & $\deg$ & $\deg$     &         & $j$    & $\deg$ & $\deg$ & $\deg$     &         \\
\hline
1 & 1 & 18 & 67.76 & 100.66 & 7.11   & 9.00 & 18 & 67.76 & 100.66 & 7.11   & 9.00 \\
  & 2 & 5  & 59.99 & 100.66 & 220.13 & 8.00 & 5  & 59.99 & 100.66 & 220.13 & 8.00 \\
2 & 1 & 3  & 40.99 & 247.27 & 277.98 & 7.00 & 3  & 40.99 & 247.27 & 277.98 & 7.00 \\
  & 2 & 19 & 44.81 & 242.00 & 201.53 & 7.00 & 12 & 44.81 & 242.00 & 141.53 & 7.00 \\
\hline \hline
    \end{tabular}
    }
    \label{tab:TC_varx}
\end{table}

Overall, the case of monitoring hurricanes via nadir-directional LEO observation satellites further exemplifies the performance of the \SDDiP algorithm relative to a realistic and highly dynamic set of hurricane scenarios. The \SDDiP algorithm converges to the \MCRP deterministic upper bound in every scenario, additionally outperforming all other \SMCRP solution methods. While the orbital maneuver sequence of the \SDDiP algorithm may not be identical to that of the \MCRP, as in the case of Scenario $18$, the sum of rewards matches precisely, indicating the existence of multiple globally optimal orbital maneuver sequences. This computational experiment not only further demonstrates the effectiveness of the \SDDiP algorithm but additionally provides a high level of performance under realistic orbital maneuvering and dynamic natural disaster monitoring parameters. Due to the slight reduction in overall problem size provided by only using two stages of reconfiguration, the difference in tradeoffs between the \SDDiP algorithm and the MDP solution methods of the \SMCRP is lessened. For instance, \textsf{VI} converges in slightly over one hour, while \textsf{QL} converges in slightly less than four hours, and the \SDDiP algorithm achieves convergence in a total of two hours. As such, the gap in both performance and computational burden is lessened due to the problem geometry of TC scenarios. However, the solution quality of the MDP solution methods may still improve through the tuning of the hyperparameter values.

\section{Conclusions} \label{sec:Conclusion}

This paper extends the Multistage Constellation Reconfiguration Problem (\MCRP) into the Stochastic MCRP (\SMCRP), incorporating stochastic target properties through the use of Stochastic Dual Dynamic Integer Programming (SDDiP) and developing the \SDDiP algorithm alongside additional solution methods. The \SDDiP algorithm presented determines the optimal policy of orbital maneuvers within an array of stochastic scenarios with realizations of uncertain mission parameters, including target visibility and rewards. Due to the nature of SDDiP, the \SDDiP algorithm converges to global optimality within finite iterations, certifying that the orbital maneuver sequence of each scenario obtains the maximum observation rewards. The \SMCRP is additionally modeled as a Markov decision process, which is then solved through the use of value iteration and Q-Learning for comparative purposes within two computational experiments.

The computational experiments conducted in Sec.~\ref{sec:Experimentation} demonstrate the capability of the \SDDiP algorithm to obtain optimal orbital maneuvering sequences under various uncertain parameter conditions for different mission purposes. Within the experiment concerning orbital targets and observation satellites in LEO, the policy obtained by the \SDDiP algorithm outperforms all other \SMCRP solution methods, especially relative to the non-maneuverable constellation and random maneuver policy. Additionally, the \SDDiP algorithm policy performance closely mirrors the deterministic upper bound \MCRP solution in a majority of scenarios. Similarly, within the experiment concerning hurricanes, the policy obtained by the \SDDiP algorithm obtains identical rewards to the \MCRP solution in all scenarios under realistically simulated random hurricane trajectories. Both experiments demonstrate the performance of the \SDDiP algorithm above other \SMCRP solution methods under various conditions, obeying orbital maneuver feasibility and physical VTW constraints while accounting for uncertain target motion and rewards. Overall, the \SDDiP algorithm provides robust orbital maneuver decisions under uncertainty, presenting itself as a valuable stochastic algorithm while leveraging the multistage reconfiguration methodology of the \MCRP. 

Future efforts to improve the optimization of orbital maneuverability under uncertain target conditions may progress in a myriad of directions. A priority would be to devise a new reformulation of the \SMCRP such that the overall mission budget parameter $c^k_{\max}$ from constraints~\eqref{MCRP:cost} may be included in both the \SDDiP algorithm and MDP formulations for a solution that may take any general orbital slots as input. Similarly, larger test instances composed of more scenarios and more targets in each scenario would further exemplify the effectiveness of the \SDDiP algorithm under a much larger set of uncertainty. Finally, an extension of the \MCRP for the EOSSP is the REOSSP in Ref.~\cite{Pearl2026REOSSP}, so a logical extension of the \SMCRP would be to incorporate other aspects of satellite scheduling also under target uncertainty, performing full operations of data downlink to ground station and solar charging with stochasticity.

\section*{Appendix A: Construction of Orbital Slots to Obey $c^k_{\max}$}

As a result of Assumption $6$ of the \SDDiP algorithm regarding stagewise independence, a special consideration is made such that the total cost of any orbital maneuver sequence does not exceed a desired budget. Specifically, the mathematical condition is that for any feasible values of $x^{sk}_{ij}$ that obey constraints~\eqref{MCRP:flow} or~\eqref{SMCRP:Flow}, the following is true: $\sum_{s \in \mathcal{S}\setminus\{0\}}\sum_{i \in \mathcal{J}^{s-1,k}} \sum_{j \in \mathcal{J}^{sk}} c^{sk}_{ij}x^{sk}_{ij} \le c^k_{\max}, ~ \forall k \in \mathcal{K}$. As such, this paper assumes the following construction method for orbital slots. Firstly, it is assumed that the same orbital slots are assigned to each stage of reconfiguration, such that $\mathcal{J}^{1k} = \mathcal{J}^{2k} = \cdots = \mathcal{J}^{Sk}$ (with the only exception being the initial condition orbital slots in $\mathcal{J}^{0k}$ being a singleton set) for each satellite $k\in\mathcal{K}$. Next, the above condition may be simplified to $S\max\left(c^{sk}_{ij}\right) \le c^k_{\max}, \forall k \in \mathcal{K}$ as the maximum cost from orbital slot $i\in\mathcal{J}^{s-1,k}$ to orbital slot $j\in\mathcal{J}^{sk}$ is the same in all $s\in\mathcal{S}\setminus\{0\}$, allowing the condition to represent the worst case cost orbital maneuver being selected in all stages of reconfiguration. Given this simplified condition, Algorithm~\ref{alg:Pruning} iteratively prunes orbital slots from a given set of orbital slots $\mathcal{J}^{sk}$ until the condition is met. The algorithm evaluates the condition before removing the offending orbital slot $j$ from the set of orbital slots in Line~5 and reducing the size of $c^{sk}_{ij}$ to match the new set of orbital slots in Line~6. Upon the removal of all orbital slots that result in a maximum cost maneuver exceeding the budget of each satellite, the new set of orbital slots and reduced cost parameters are returned. 

\begin{algorithm}[!ht]
    \DontPrintSemicolon
    \caption{Pruning orbital slots to obey a given budget.}
    \label{alg:Pruning}
    \textbf{Input:} $\mathcal{J}^{sk}, ~ c^{sk}_{ij}, ~ c^k_{\max}$ \;
    \For{$k = 1, \ldots, K$}{
        \While{$S\max\left(c^{sk}_{ij}\right) > c^k_{\max}$}{
            $[i, j] = \underset{i,j}{\arg\max} \left(c^{sk}_{ij}\right)$\;
            $\mathcal{J}^{sk} \gets \mathcal{J}^{sk} \setminus\{j\}$ \;
            Recompute $c^{sk}_{ij}$ according to updated $\mathcal{J}^{sk}$ \;
        }
    }
    \Return $\mathcal{J}^{sk}, ~c^{sk}_{ij}$
\end{algorithm}

This paper considers a different set of orbital slots in each computational experiment in Sec.~\ref{sec:Experimentation}, one of which concerns orbital slots that vary only in the argument of latitude, while the other considers additional changes in inclination and RAAN. As mentioned in Sec.~\ref{sec:Experimentation}, a total of seven orbital maneuvers are possible, for which the cost $c^{sk}_{ij}$ is computed through the use of the circular coplanar phasing problem and other associated one- or two-impulse direct transfer analytical algorithms from Ref.~\cite{Vallado2022}. As a result, $c^{sk}_{ij} = 0$ when $i = j$, avoiding an edge case in Algorithm~\ref{alg:Pruning} where all costs exceed the budget and $\mathcal{J}^{sk} = \emptyset$. Therefore, given a budget of $c^k_{\max} = \SI{2.5}{\kilo\meter}, ~ \forall k \in \mathcal{K}$, the first experiment using the argument of latitude orbital slots (shown in Fig.~\ref{fig:orbital_targets}) has a maximum cost of $\max\left( c^{sk}_{ij} \right) = 0.62$ with the condition $S\left(0.62\right) = 2.49 < c^k_{\max},~\forall k \in \mathcal{K}$. Similarly, the second experiment with orbital slots additionally varying in inclination and RAAN (shown in Fig.~\ref{fig:PlanarSlots}) has a maximum cost of $\max\left( c^{sk}_{ij} \right) = 1.09$ with the condition $S\left(1.09\right) = 2.18 < c^k_{\max},~\forall k \in \mathcal{K}$. As such, the orbital slots employed in the computational experiments obey the assumption that no sequence of maneuvers exceeds the established budget, satisfying stagewise independence and the maximum maneuver budget. 

\bibliography{References_V3}

\end{document}